\documentclass[11pt,letterpaper]{article}

\usepackage{customize}

\title{Interpolating Log-Determinant and Trace of the Powers of Matrix \(\tens{A} + \t \tens{B}\)}

\author{Siavash Ameli\protect\thanks{Email address: \href{mailto:sameli@berkeley.edu}{\protect\nolinkurl{sameli@berkeley.edu}}}~}
\author{Shawn C. Shadden\protect\thanks{Email address: \href{mailto:shadden@berkeley.edu}{\protect\nolinkurl{shadden@berkeley.edu}}}}

\affil{\small\textit{Mechanical Engineering}, \small\textit{University of California}, \small\textit{Berkeley, CA, USA 94720}}

\date{}

\def \eig {\lambda}

\def \t {t}

\DeclareMathOperator*{\argmin}{arg\,min}

\ifpdf
\hypersetup{
  pdftitle={Interpolating Log-Determinant and Trace of the Powers of Matrix A+tB},
  pdfauthor={S. Ameli and S. C. Shadden}
}
\fi

\def\CC{{C\nolinebreak[4]\hspace{-.05em}\raisebox{.4ex}{\tiny\bf ++}}}

\DeclareFixedFont{\ttb}{T1}{txtt}{m}{n}{9} 
\DeclareFixedFont{\ttm}{T1}{txtt}{m}{n}{9}  

\definecolor{deepblue}{rgb}{0, 0, 0.5}
\definecolor{deepred}{rgb}{0.6, 0, 0}
\definecolor{deepgreen}{rgb}{0.4, 0.4, 0.4}
\definecolor{cmnt}{HTML}{778899}

\lstdefinestyle{myStyle}{
    language=Python,
    basicstyle=\small\ttfamily,
    morekeywords={self},              
    keywordstyle=\small\color{deepblue},
    emph={InterpolateTraceinv, imate, interpolate},          
    emphstyle=\small\ttfamily\bfseries\color{deepred},    
    commentstyle=\small\color{cmnt}\ttfamily\itshape,
    stringstyle=\color{deepgreen},
    frame=tb,                         
    showstringspaces=false,
    escapeinside={;\#}{\#;},           
    escapebegin=\color{cmnt},        
    numberstyle = \tiny\color{black},
    numbers=left,
    firstnumber=1,
    stepnumber=5,
    columns=flexible,
}

\begin{document}

\maketitle

\begin{abstract}
    We develop heuristic interpolation methods for the functions \(\t \mapsto \log \det \left( \tens{A} + \t \tens{B} \right)\) and \(\t \mapsto \trace\left( (\tens{A} + \t \tens{B})^{p} \right)\) where the matrices \(\tens{A}\) and \(\tens{B}\) are Hermitian and positive (semi) definite and \(p\) and \(\t\) are real variables. These functions are featured in many applications in statistics, machine learning, and computational physics. The presented interpolation functions are based on the modification of sharp bounds for these functions. We demonstrate the accuracy and performance of the proposed method with numerical examples, namely, the marginal maximum likelihood estimation for Gaussian process regression and the estimation of the regularization parameter of ridge regression with the generalized cross-validation method.


\end{abstract}


\section{Introduction}


Estimation of the determinant and trace of matrices is a key component and often a computational challenge in many algorithms in data analysis, statistics, machine learning, computational physics, and computational biology. Some applications of trace estimation can be found in \citet{UBARU-2018}. A few examples of such applications are high-performance uncertainty quantification \citep{BEKAS-2012,KALANTZIS-2013}, optimal Bayesian experimental design \citep{CHALONER-1995}, regression using Gaussian process \citep{MACKAY-2003}, rank estimation \citep{UBARU-2016}, and computing observables in lattice quantum chromodynamics \citep{WU-2016}. 


\subsection{Motivation} \label{sec:Motivation}

In this paper, we are interested in estimating the functions
\begin{subequations}
\begin{align}
    &\t \mapsto \log \det \left( \tens{A} + \t \tens{B} \right), \label{eq:det-function}
    \intertext{and}
    &\t \mapsto \trace \left( (\tens{A} + \t \tens{B})^{p} \right), \label{eq:trace-function}
\end{align}
\end{subequations}
where \(\tens{A}\) and \(\tens{B}\) are Hermitian and positive semi-definite (positive-definite if \(p < 0\)), and \(p\) and \(\t\) are real numbers\footnote{We use boldface lowercase letters for vectors, boldface upper case letters for matrices, and normal face letters for scalars, including the components of vectors and matrices, such as \(x_i\) and \(H_{ij}\) respectively for the components of the vector \(\vect{x}\) and the matrix \(\tens{H}\).}. These functions are featured in a vast number of applications in statistics and machine learning. Often, in these applications, the goal is to optimize a problem for the parameter \(\t\), and the above functions should be evaluated for a wide range of \(\t\) during the optimization process.

A common example of such an application can be found in regularization techniques applied to inverse problems and supervised learning. For instance, in ridge regression by generalized cross-validation \citep{WAHBA-1977,CRAVEN-1978,GOLUB-1997}, the optimal regularization parameter \(\t\) is sought by minimizing a function that involves \eqref{eq:trace-function} at \(p = -1\) (see \Cref{sec:example-2}). Another common usage of \eqref{eq:det-function} and \eqref{eq:trace-function}, for instance, is the mixed covariance functions of the form \(\tens{A} + \t \tens{I}\) that appear frequently in Gaussian processes with additive noise \citep{AMELI-2022-a, AMELI-2022-d} (see also \Cref{sec:example-1}). In most of these applications, the log-determinant of the covariance matrix is common, particularly in likelihood functions or related variants. Namely, if one aims to maximize the likelihood by its derivative with respect to the parameter, the expression,
\begin{equation*}
    \frac{\partial}{\partial \t} \log \det (\tens{A} + \t \tens{I}) = \trace \left( (\tens{A} + \t \tens{I})^{-1} \right),
\end{equation*}
frequently appears. More generally, the function \eqref{eq:trace-function} for \(p \in \mathbb{Z}_{< 0}\) appears in the \(\vert p \vert\)-th derivative of such likelihood functions. Other examples of \eqref{eq:det-function} and \eqref{eq:trace-function} are in experimental design \citep{HABER-2008}, probabilistic principal component analysis \citep[Sec. 12.2]{BISHOP-2006}, relevance vector machines \citep{TIPPING-2001} and \citep[Sec. 7.2]{BISHOP-2006}, kernel smoothing \citep[Sec. 2.6]{RASMUSSEN-2006}, and Bayesian linear models \citep[Sec. 3.3]{BISHOP-2006}.


\subsection{Overview of Related Works} \label{sec:related-works}

The difficulty of estimating \eqref{eq:det-function} and \eqref{eq:trace-function} in all the above applications is that the matrices are generally large. Also, often in these applications, cases of particular interest in \eqref{eq:trace-function} are when \(p < 0\), but the \(\vert p \vert\)-th inverse of the matrix \(\tens{A} + t \tens{B}\) is not available explicitly, rather it is implicitly known by matrix-vector multiplications through solving a linear system. Because of these, the evaluation of \eqref{eq:det-function}  and \eqref{eq:trace-function} are usually the main computational challenge in these problems, and several algorithms have been developed to address this problem.

The determinant and trace of the inverse of a Hermitian and positive-definite matrix can be calculated by the Cholesky factorization (cf.\ equation \eqref{eq:det-chol} and \eqref{eq:trace-chol} in \Cref{sec:example-1}). Using the Cholesky factorization, \citet{TAKAHASHI-1973} developed a method to find desired entries of a matrix inverse, such as its diagonals. The latter method was extended by \citet{NIESSNER-1983}. Also, \citet{GOLUB-1980} found entries of the inverse of the covariance matrix provided that the corresponding entries of its Cholesky factorization are non-zero. The complexity of this method is \(\mathcal{O}(nw)\) where \(w\) is the bandwidth of the Cholesky matrix (see also \citet[Sec. 6.7.4]{BJORK-1996}). Recently, probing and hierarchical probing methods were presented by \citet{TANG-2012} and \citet{STATHOPOULOS-2013}, respectively, to compute the diagonal entries of a matrix inverse.

In contrast to the above exact methods, many approximation methods have been developed. The stochastic trace estimator by \citet{HUTCHINSON-1990}, which evolved from \citet{GIRARD-1989}, uses Monte-Carlo sampling of random vectors with a Gaussian or Rademacher distribution. A similar concept was presented by \citet{GIBBS-1997}. Another randomized trace estimator was given by \citet{AVRON-2011} for symmetric and positive-definite implicit matrices. Based on the stochastic trace estimation, \citet{WU-2016} interpolated the diagonals of a matrix inverse. Also, \citet{SAIBABA-2017} improved the randomized estimation by a low-rank approximation of the matrix. Another tier of methods combines the idea of a stochastic trace estimator and Lanczos quadrature \citep{GOLUB-1994,BAI-1996,BAI-1997,GOLUB-2009}, which is known as stochastic Lanczos quadrature (SLQ). The numerical details of the SLQ method using either Lanczos tridiagonalization or Golub-Kahn bidiagonalization can be found for instance in \citet[Algorithms 1 and 2]{UBARU-2017}.


\subsection{Objective and Our Contribution}

Our objective is to develop a method to efficiently estimate  \eqref{eq:det-function} or \eqref{eq:trace-function} for a wide range of \(\t\). Note, if \(\tens{B}\) is the identity matrix and the matrix \(\tens{A}\) is small enough to pre-compute is eigenvalues, \(\eig_i(\tens{A})\), then, the evaluation of \eqref{eq:det-function} and \eqref{eq:trace-function} for any \(\t\) is immediate by
\begin{subequations}
\begin{align}
    &\log \det (\tens{A} + \t \tens{I}) = \sum_{i = 1}^n \log(\eig_i(\tens{A}) + \t), \label{eq:det-AI} \\
    &\trace\left( (\tens{A} + \t \tens{I})^{p} \right) = \sum_{i = 1}^n (\eig_i(\tens{A}) + \t)^{p}. \label{eq:tr-AI}
\end{align}
\end{subequations}
However, for large matrices, estimating all eigenvalues is impractical. Thus, we herein develop an interpolation scheme for the functions \eqref{eq:det-function} and \eqref{eq:trace-function} based on the following developments:

\begin{itemize}
    \item We present a Schatten-type operator that unifies the representation of \eqref{eq:det-function} and \eqref{eq:trace-function} by a single continuous function. This operator leads to definitions of an associated norm and anti-norm on matrices. Sharp inequalities for this norm and anti-norm on the sum of two Hermitian and positive (semi) definite matrices provide rough estimates for \eqref{eq:det-function} and \eqref{eq:trace-function}.

    \item We propose two interpolation methods based on the sharp norm and anti-norm inequalities mentioned above. Namely, we introduce interpolation functions based on a linear combination of orthogonal basis functions, or interpolation by rational polynomials.
\end{itemize}

We demonstrate the computational advantage of our method through two examples:

\begin{itemize}
    \item \emph{Gaussian process regression.} We compute \eqref{eq:det-function} and \eqref{eq:trace-function} in the context of marginal likelihood estimation of Gaussian process regression. We show that with very few interpolation points, an accuracy of \(0.01 \% \) can be achieved.

    \item \emph{Ridge regression.} We estimate the regularization parameter of ridge regression with the generalized cross-validation method. We demonstrate that with only a few interpolation points, the ridge parameters can be estimated and the overall computational cost is reduced by 2 orders of magnitude.
\end{itemize}

The outline of the paper is as follows. In \Cref{sec:ineq}, we present matrix determinant and trace inequalities. In \Cref{sec:interpolation}, we propose interpolation methods. In \Cref{sec:example} we provide examples and a software package that implements the presented algorithms. In \Cref{sec:other-app}, we provide further applications of the method. \Cref{sec:conclusion} concludes the paper. Proofs are given in \Cref{sec:pf}.



\section{Determinant and Trace Inequalities} \label{sec:ineq}

We will derive interpolations for \eqref{eq:det-function} and \eqref{eq:trace-function} by modifying sharp bounds for these functions. In this section, we present these bounds. Without loss of generality, we temporarily omit the parameter \(\t\). However, in \Cref{sec:interpolation}, we will retrieve the desired relations by replacing \(\tens{B}\) with \(\vert \t \vert \tens{B}\).

Let \(\mathcal{M}_{n, m}(\mathbb{C})\) denote the space of all \(n \times m\) matrices with entries over the field \(\mathbb{C}\). We assume \(\tens{A}, \tens{B} \in \mathcal{M}_{n, n}(\mathbb{C})\) are Hermitian and positive semi-definite. Furthermore, for \(p < 0\), we require matrices \(\tens{A}\) and \(\tens{B}\) to be positive-definite. The notations \(\tens{A} \succ \tens{B}\) and \(\tens{A} \succeq \tens{B}\) on matrices \(\tens{A}\) and \(\tens{B}\) denotes \(\tens{A} - \tens{B}\) is positive-definite and positive semi-definite, respectively. Also, \(\eig(\tens{A}) \coloneqq (\eig_1(\tens{A}), \cdots, \eig_n(\tens{A}))\) indicates the \(n\)-tuple of eigenvalues of matrix \(\tens{A}\).

Define a Schatten-class operator \(\| \cdot \|_p: \mathcal{M}_{n, n}(\mathbb{C}) \mapsto \mathbb{R}_{\geq0}\) by
\begin{equation}
    \| \tens{A} \|_p \coloneqq
    \begin{cases}
        \left( \det(\vert \tens{A} \vert) \right)^{\frac{1}{n}}, & p = 0, \\
        \left( \frac{1}{n} \trace(\vert \tens{A}\vert^p) \right)^{\frac{1}{p}}, & p \in \mathbb{R} \setminus \{0\},
    \end{cases}
        \label{eq:schatten}
\end{equation}
where \(\vert \tens{A} \vert \coloneqq \sqrt{\tens{A}^* \tens{A}}\) and \(\tens{A}^*\) denotes the Hermitian transpose of \(\tens{A}\). Since we assume the matrices are Hermitian and at least positive semi-definite, we omit \(\vert \cdot \vert\) in subsequent expressions. Also, we note that \(p \mapsto \| \cdot \|_p\) is continuous at \(p = 0\) since
\begin{equation}
    \| \tens{A} \|_0 = \lim_{p \to 0} \| \tens{A} \|_p. \label{eq:lim}
\end{equation}
Namely, \eqref{eq:lim} is justified by observing that \(\| \tens{A} \|_p = M_p(\eig(\tens{A})) \), where \(M_p(\eig(\tens{A}))\) is the generalized mean of \(\eig(\tens{A})\) defined by
\begin{equation}
    M_p(\eig(\tens{A})) \coloneqq
    \begin{cases}
        \Big( \prod_{i=1}^n \eig_i(\tens{A}) \Big)^{\frac{1}{n}}, & p = 0, \\
        \Big( \frac{1}{n} \sum_{i=1}^n \eig_i^p(\tens{A}) \Big)^{\frac{1}{p}}, & p \in \mathbb{R} \setminus \{0\}.
    \end{cases} \label{eq:gen-mean}
\end{equation}
It is known that the generalized mean converges to the geometric mean, \(M_0\), as \(p \to 0\) \citep[p. 15]{HARDY-1952}, which concludes \eqref{eq:lim}. 

For \(p \in [1, \infty)\), the operator \(\| \cdot \|_p\) is an equivalent norm to the Schatten \(p\)-norm of \(\tens{A}\). Conventionally, the Schatten norm is defined without the normalizing factor \(\frac{1}{n}\) in \eqref{eq:schatten}, but the inclusion of this factor is justified by the continuity granted by \eqref{eq:lim}. The Schatten norm is a subadditive function, meaning that it satisfies the triangle inequality
\begin{subequations}
\begin{equation}
    \| \tens{A} + \tens{B} \|_{p} \leq \| \tens{A} \|_{p} + \| \tens{B} \|_{p}. \label{eq:sch-ineq}
\end{equation}
The reverse triangle inequality follows from the above by
\begin{equation}
    \| \tens{A} - \tens{B} \|_{p} \geq \| \tens{A} \|_{p} - \| \tens{B} \|_{p}, \label{eq:sch-ineq-rev}
\end{equation}
\end{subequations}
provided that \(\tens{A} \succeq \tens{B}\) for \eqref{eq:sch-ineq-rev} to hold. For \(p = 1\), the two above relations become equality by the additivity of the trace operator.

When \(p < 1\), the operator \eqref{eq:schatten} is no longer a norm, rather, is an anti-norm \citep{BOURIN-2011-a} that satisfies superadditivity property
\begin{subequations}
\begin{equation}
    \| \tens{A} + \tens{B} \|_{p} \geq \| \tens{A} \|_{p} + \| \tens{B} \|_{p}. \label{eq:anti-ineq}
\end{equation}
A reversed inequality can also be derived from the above as
\begin{equation}
    \| \tens{A} - \tens{B} \|_{p} \leq \| \tens{A} \|_{p} - \| \tens{B} \|_{p}, \label{eq:anti-ineq-rev}
\end{equation}
\end{subequations}
provided that \(\tens{A} \succeq \tens{B}\) (or \(\tens{A} \succ \tens{B}\) if \(p < 0\)) for \eqref{eq:anti-ineq-rev} to hold. We also note that inequality \eqref{eq:anti-ineq} at \(p=0\) reduces to Brunn-Minkowski determinant inequality \citep[p. 482, Theorem 7.8.8]{HORN-1990}
\begin{equation}
    \left(\det (\tens{A} + \tens{B}) \right)^{\frac{1}{n}} \geq \left(\det (\tens{A}) \right)^{\frac{1}{n}} + \left(\det (\tens{B}) \right)^{\frac{1}{n}}.
\end{equation}

For proofs and discussions of a general class of anti-norms, which includes \eqref{eq:schatten}, we refer the reader to \citep{BOURIN-2011-a} and \citep{BOURIN-2014}. However, in \Cref{sec:pf}, we provide a direct proof of \eqref{eq:anti-ineq} and \eqref{eq:anti-ineq-rev} and the necessary and sufficient conditions for equality to hold in these relations for the operator \eqref{eq:schatten} at \(p < 1\). 

\begin{remark}[Comparisons to other inequalities] \label{rem:lower-bound}
    There are other known bounds to functions \eqref{eq:det-function} and \eqref{eq:trace-function}. For instance, for the common case of \(p=-1\), we can obtain the upper bound \cite[p. 210, Theorem 7.7]{ZHANG-2011} 
    \begin{equation}
        \trace\left( (\tens{A} + \tens{B})^{-1} \right) \leq 
        \frac{1}{4} \left( \trace(\tens{A}^{-1}) + \trace(\tens{B}^{-1}) \right). \label{eq:trace-upper-bound}
    \end{equation}
    Also, a lower bound can be obtained, for instance, by the arithmetic-harmonic mean inequality \(M_{-1}(\eig(\tens{A} \pm \tens{B})) \leq M_1(\eig(\tens{A} \pm \tens{B}))\), where \(M_{-1}\) and \(M_1\) are the harmonic mean and arithmetic mean, respectively, \cite[Ch. 2, Theorem 1]{MITRINOVIC-1970}, which leads to
    \begin{equation}
        \trace((\tens{A} \pm \tens{B})^{-1}) \geq \frac{n^2}{\trace(\tens{A}) \pm \trace(\tens{B})}. \label{eq:trace-lower-bound}
    \end{equation}
    The inequalities \eqref{eq:trace-upper-bound} or \eqref{eq:trace-lower-bound}, however, are not as useful as the inequality \eqref{eq:anti-ineq} for \(p=-1\), since if \(\tens{B}\) is either too small or too large compared to \(\tens{A}\), \eqref{eq:trace-upper-bound} and \eqref{eq:trace-lower-bound} do not asymptote to equality, whereas \eqref{eq:anti-ineq} and \eqref{eq:anti-ineq-rev} become asymptotic equalities, which is a desired property for our purpose.
\end{remark}


\section{Interpolation of Determinant and Trace} \label{sec:interpolation}

We use the bounds provided by inequalities \eqref{eq:sch-ineq}, \eqref{eq:sch-ineq-rev}, \eqref{eq:anti-ineq}, and \eqref{eq:anti-ineq-rev} to interpolate the functions \eqref{eq:trace-function} and \eqref{eq:det-function}. To this end, we replace the matrix \(\tens{B}\) with \(\vert \t \vert \tens{B}\) in the bounds found in \Cref{sec:ineq}. Define
\begin{equation*}
    \tau_p(\t) \coloneqq \frac{\| \tens{A} + \t \tens{B} \|_p}{\| \tens{B}\|_p},
    \qquad
    \text{and}
    \qquad
    \tau_{p, 0} \coloneqq \tau_p(0).
\end{equation*}
We assume \(\tau_{p, 0}\) is known by directly computing \((\frac{1}{n}\trace(\tens{A}^p))^{\frac{1}{p}}\) and \((\frac{1}{n}\trace(\tens{B}^p))^{\frac{1}{p}}\) when \(p \neq 0\), or \((\det(\tens{A}))^{\frac{1}{n}}\) and \((\det(\tens{B}))^{\frac{1}{n}}\) when \(p = 0\)\footnote{Computing the determinant directly should be avoided as it can be a very large number. Rather, \((\det(\cdot))^{\frac{1}{n}}\) can be computed via \(\exp(\frac{1}{n}\log \det (\cdot))\). See \Cref{sec:example-1} for an example.}. Then, \eqref{eq:sch-ineq} and \eqref{eq:anti-ineq} imply
\begin{subequations}
\begin{alignat}{3}
    &\tau_p(\t) \leq \tau_{p, 0} + \t, &\quad p \geq 1, &\quad \t \in [0,\infty), \label{eq:tau-ineq-1} \\
    &\tau_p(\t) \geq \tau_{p, 0} + \t, &\quad p < 1,    &\quad \t \in [0,\infty), \label{eq:tau-ineq-2} \\
    \intertext{and \eqref{eq:sch-ineq-rev} and \eqref{eq:anti-ineq-rev} imply}
    &\tau_p(\t) \geq \tau_{p, 0} + \t, &\quad p \geq 1, &\quad \t \in (\t_{\inf},0], \label{eq:tau-ineq-neg-1} \\
    &\tau_p(\t) \leq \tau_{p, 0} + \t, &\quad p < 1,    &\quad \t \in (\t_{\inf},0],
    \label{eq:tau-ineq-neg-2}
\end{alignat}
\end{subequations}
where \(t_{\inf} \coloneqq \inf \{t ~\vert \tens{A} + t \tens{B} \succ \tens{0}\}\). The above sharp inequalities become equality at \(\t = 0\). Also, \eqref{eq:tau-ineq-1} and \eqref{eq:tau-ineq-2} become asymptotic equalities as \(\t \to \infty\). Based on the above, the bound function
\begin{equation}
    \hat{\tau}_p(\t) \coloneqq \tau_{p,0} + \t, \label{eq:tau-upper}
\end{equation}
can be regarded as a reasonable approximation of \(\tau_p(\t)\) at \(\vert \t \vert \ll \tau_{p, 0}\) where \(\tau_p(\t) \approx \tau_{p, 0}\), and at \(\t \gg \tau_{p, 0}\) where \(\tau_p(\t) \approx \t\). We expect \(\hat{\tau}_p(\t)\) to deviate the most from \(\tau_p(\t)\) when \(\mathcal{O}(\t \tau_{p, 0}^{-1}) \approx 1\).

Furthermore, to improve the approximation in the intermediate interval \(\t \tau_{p, 0}^{-1} \in (c,c^{-1})\) for some \(c \ll 1\), we define interpolating functions based on the above bounds to honor the known function values at some intermediate points \(\t_i \in (c \tau_{p, 0},c^{-1} \tau_{p, 0})\). In particular, we specify interpolation points over logarithmically spaced intervals, because \(\t\) is usually varied in a wide logarithmic range in most applications. We compute the function values at the interpolation points, \(\tau_p(\t_i)\), with any of the trace estimation methods mentioned earlier.

Many types of interpolating functions can be employed to improve the above approximation. However, we seek interpolating functions whose parameters can be easily obtained by solving a linear system of equations. We define two such types of interpolations, namely, by a linear combination of basis functions and by rational polynomials, respectively in \Cref{sec:interp-1} and \Cref{sec:interp-2}.


\subsection{Interpolation with a Linear Combination of Basis Functions} \label{sec:interp-1}

Based on \eqref{eq:tau-upper}, we define an interpolating function \(\tilde{\tau}_p(\t)\) by
\begin{equation}
    \tilde{\tau}_p(\t) \coloneqq \tau_{p, 0} + \sum_{i = 0}^{q} w_i \phi_i(\t), \label{eq:approx}
\end{equation}
where \(\phi_i\) are basis functions and \(w_i\) the weights. The basis functions
\begin{equation}
    \phi_i(\t) = \t^{\frac{1}{i+1}}, \qquad i = 0,\dots,q, \label{eq:basis}
\end{equation}
for the domain \(\t \in [0, \infty)\) can be used, which are inverse functions of the monomials and we refer to them as inverse-monomials. These basis functions satisfy the conditions  \(\phi_0(\t) = \t\), \(\phi_i(0) = 0\), and \(\phi_0(\t) \gg \phi_{i}(\t)\), \(i > 0\) when \(\t \gg 1\). For consistency with \eqref{eq:tau-upper}, we set \(w_0 = 1\).  The coefficients \(w_i\), \(i = 1,\dots,q\) are found by solving a linear system of \(q\) equations using a priori known values \(\tau_{p, i} \coloneqq \tau_p(\t_i)\), \(i = 1,\dots,q\). When \(q = 0\), no intermediate interpolation point is introduced and the approximation function is the same as the bound \(\hat{\tau}_p(\t)\) given by \eqref{eq:tau-upper}.

\begin{remark}
    An alternative could be to use monomials \(\t^i\) for interpolation functions, \eg
    \begin{equation}
        \tilde{\tau}_p(\t)^{q+1} \coloneqq \tau_{p, 0}^{q+1} + \sum_{i = 1}^{q+1} w_i \t^i, \label{eq:approx-2}
    \end{equation}
    with \(w_{q+1} = 1\), and the rest of the weights \(w_i\), \(i = 1,\dots,q\) determined from the known values of the function. 
    This is not particularly useful in practice, as the exponentiation terms, \(\t^i\), cause arithmetic underflows; also, Runge's phenomenon occurs even for low-order interpolations  \(q > 1\).
\end{remark}


%

In practice, just a few interpolating points \(\t_i\) are sufficient to obtain a reasonable interpolation of \(\tau_p(\t)\). However, when more interpolation points are used (such as when \(p \geq 6\)), the linear system of equations for the weights \(w_i\) becomes ill-conditioned. To overcome this issue, orthogonal basis functions can be used (see \eg \citet[Sec. 7.1]{SEBER-2012} for a general discussion).

For our application, we seek basis functions \(\phi_i^{\perp}(\t)\) that are orthogonal on the unit interval \(\t \in [0,1]\). Since we are interested in functions in the logarithmic scale of \(\t\), we define the inner product in the space of functions using the Haar measure \(\mathrm{d} \log(\t) = \mathrm{d} \t / \t\) . Applicability of Haar measure can be justified by letting \(\t_i = e^{x_i}\), where \(x_i\) are normally spaced interpolant points. Following the discussion of \citet[Sec. 7.1]{SEBER-2012} for linear regression using orthogonal polynomials, we use the conventional integrals with the Lebesgue measure \(\mathrm{d}x\) to define the inner product of functions. The measure \(\mathrm{d} x\) is equivalent to the Haar measure \(\mathrm{d} \log \t\) for the variable \(\t\).

\begin{figure}[bt!]
    \centering
    \includegraphics[width=0.6\textwidth]{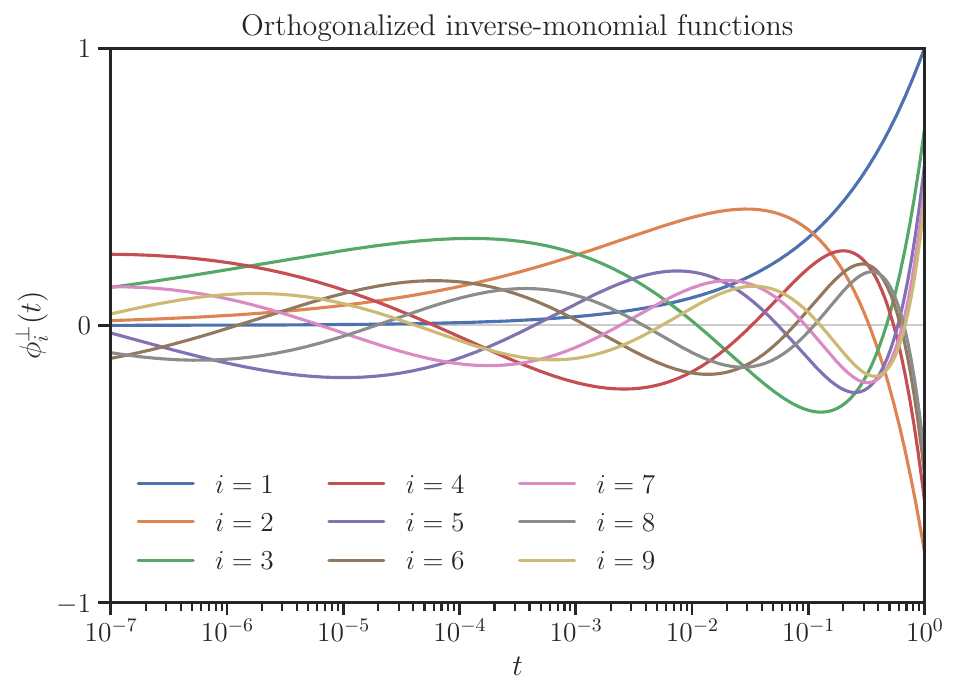}
    \caption{Orthogonalized inverse-monomial functions \(\phi_i^{\perp}(\t)\) in the logarithmic scale of \(\t\).}
    \label{fig:OrthogonalFunctions}
\end{figure}

The desired orthogonality condition in the Hilbert space of functions on \([0,1]\) with respect to the Haar measure becomes
\begin{equation}
    \langle \phi_i^{\perp}, \phi_j^{\perp} \rangle_{L^2([0,1],\mathrm{d}\t/\t)} = \int_{0}^1 \phi_i^{\perp}(\t) \phi_j^{\perp}(\t) \frac{\mathrm{d} \t}{\t} = \delta_{ij}, \label{eq:inner-prod}
\end{equation}
where \(\delta_{ij}\) is the Kronecker delta function. A set of orthogonal functions \(\phi_i^{\perp}(\t)\) can be constructed from the set of non-orthogonal basis functions \(\{\phi_i\}_{i=1}^q\) in \eqref{eq:basis} by recursive application of Gram-Schmidt orthogonalization 
\begin{equation}
    \phi_i^{\perp}(\t) = \alpha_i \sum_{j = 1}^q a_{ij} \phi_j(\t), \qquad i = 1,\dots,q. \label{eq:ortho-funct}
\end{equation}
The first nine orthogonal basis functions are shown in \Cref{fig:OrthogonalFunctions} and the respective coefficients \(\alpha_i\) and \(a_{ij}\) are given by \Cref{table:ortho-funct}.\footnote{We developed the python package \texttt{ortho} to generate an arbitrary number of orthogonal functions \(\phi_j^{\perp}(\t)\) using symbolic computations. See \url{https://ameli.github.io/ortho} for details.}

\begin{table}[hbt!]
    \centering
    \caption{Coefficients of orthogonal functions in \eqref{eq:ortho-funct}} \label{table:ortho-funct}
    \begin{adjustbox}{width=1\textwidth}
    \begin{small}
    \begin{tabular}{c| c c c c c c c c c c}
        \toprule
        $i$ & $\alpha_i$ & $a_{i1}$ & $a_{i2}$ & $a_{i3}$ & $a_{i4}$ & $a_{i5}$ & $a_{i6}$ & $a_{i7}$ & $a_{i8}$ & $a_{i9}$ \\
        \midrule
        $1$ & $+\sqrt{2/2}$  & 1     &          &         &           &          &            &          &           &         \\
        $2$ & $-\sqrt{2/3}$  & $6$   & $-5$     &         &           &          &            &          &           &         \\
        $3$ & $+\sqrt{2/4}$  & $20$  & $-40$    & $21$    &           &          &            &          &           &         \\
        $4$ & $-\sqrt{2/5}$  & $50$  & $-175$   & $210$   & $-84$     &          &            &          &           &         \\
        $5$ & $+\sqrt{2/6}$  & $105$ & $-560$   & $1134$  & $-1008$   & $330$    &            &          &           &         \\
        $6$ & $-\sqrt{2/7}$  & $196$ & $-1470$  & $4410$  & $-6468$   & $4620$   & $-1287$    &          &           &         \\
        $7$ & $+\sqrt{2/8}$  & $336$ & $-3360$  & $13860$ & $-29568$  & $34320$  & $-20592$   & $5005$   &           &         \\
        $8$ & $-\sqrt{2/9}$  & $540$ & $-6930$  & $37422$ & $-108108$ & $180180$ & $-173745$  & $90090$  & $-19448$  &         \\
        $9$ & $+\sqrt{2/10}$ & $825$ & $-13200$ & $90090$ & $-336336$ & $750750$ & $-1029600$ & $850850$ & $-388960$ & $75582$ \\
        \bottomrule
    \end{tabular}
    \end{small}
    \end{adjustbox}
\end{table}

A set of orthogonal functions can also be defined on intervals other than \([0,1]\) by adjusting the bounds of integration in \eqref{eq:inner-prod}, which yields a different set of function coefficients. However, it is more convenient to fix the domain of orthogonal functions to the unit interval \([0,1]\), and later scale the domain as desired, \eg to \([0,l]\) where \(l \coloneqq \operatorname{max}(\t_i)\). Although this approach does not lead to orthogonal functions in \([0,l]\), it nonetheless produces a well-conditioned system of equations for the weights \(w_i\).

\begin{remark} \label{rem:oscillate}
    The interpolation function defined in \eqref{eq:approx} asymptotes consistently to \(\tilde{\tau}_p(\t) \to \t\) at \(\t \gg \tau_{p, 0}\). On the other end, the convergence \(\tilde{\tau}_p(\t) \to \tau_{p, 0}\) at \(\t \ll \tau_{p, 0}\) is not uniform, rather the interpolation function oscillates. This behavior is originated from the basis functions \(\phi_i\), \(i > 0\), that are not independent at \(\t \ll \tau_{p, 0}\), particularly, near the origin. This dependency of basis functions cannot be resolved by the orthogonalized functions \(\phi_i^{\perp}\), as they are orthogonal with respect to the singular weight function \(\t^{-1}\) at the origin. Thus, \eqref{eq:approx} should not be employed on very small logarithmic scales, rather, other interpolation functions should be employed for such purpose, such as presented in \Cref{sec:interp-2}.
\end{remark}


\subsection{Interpolation with Rational Polynomials} \label{sec:interp-2}

We define another type of interpolating function that can perform well at small scales of \(\t\), by using rational polynomials. Define
\begin{equation}
    \tilde{\tau}_p(\t) \coloneqq \frac{\t^{q+1} + a_{q} \t^{q} + \cdots + a_1 \t + a_0}{\t^{q} + b_{q-1} \t^{q-1} + \cdots + b_1 \t + b_0}, \label{eq:rational}
\end{equation}
which is the Pad\'{e} approximation of \(\tau_p\) of order \([q+1, q]\). We set \(a_0 = b_0 \tau_{p, 0}\) in order to satisfy \(\tau_p(0) = \tau_{p, 0}\). Also, we note that the above interpolation also satisfies the asymptotic relation \(\tau_p(\t) \to \t \) as \(\t \to \infty\). At \(q = 0\), when no interpolant point is used, the above interpolation function falls back to \eqref{eq:tau-upper} by setting \(b_0 = 1\). For \(q > 0\), \(2q\) interpolant points \(\t_i\) are needed to solve the linear system of equations for the coefficients \(a_1,\dots,a_{q}\) and \(b_0,\dots,b_{q-1}\).

An alternative rational polynomial is the Chebyshev rational function \citep{GUO-2002}
\begin{equation}
    r_i(\t) \coloneqq T_i\left( \frac{\t-1}{\t+1} \right), \label{eq:cheb}
\end{equation}
where \(T_i\), \(i \in \mathbb{N}\) are the Chebyshev polynomials of the first kind. The Chebyshev rational functions are orthogonal in \([0, \infty)\) with respect to the weight function \((\t+1)^{-1} \sqrt{\t}\) and satisfy the recursive relation \(r_{i+1}(\t) = 2((\t-1)/(\t+1)) r_i(\t) - r_{i-1}(\t)\) with \(r_0(t) = 1\) and \(r_1(\t) = (\t-1)/(\t+1)\). An interpolation of \(\tau_p(t)\) using Chebyshev rational functions can be given by
\begin{equation}
    \frac{\tilde{\tau}_p(\t)}{\tau_{p, 0} + \t} - 1 = \sum_{i=1}^{q+1} \frac{w_i}{2} \left(1-r_i\left(\frac{\t}{\alpha} \right) \right), \label{eq:cheb-interp}
\end{equation}
where \(\alpha > 0\) is a given scale parameter and will be explained shortly. Both sides of the above relation converge to zero at \(\t \to \infty\). To satisfy \(\tau_{p}(0) = \tau_{p, 0}\), we require \(\sum_{i=1}^{q+1} \frac{w_i}{2} (1 - (-1)^i) = 0\) considering \(r_i(0) = (-1)^i\). The latter condition together with \(q\) linear equations on the interpolating points \(t_i > 0\), \(i=1, \dots, q\) solve the weights \(w_i\). An advantage of using the above interpolation scheme is that we can arrange the interpolant points \(\t_i\) on the corresponding Chebyshev nodes to reduce the interpolation error.

\begin{figure}[bt!]
    \centering
    \includegraphics[width=0.6\textwidth]{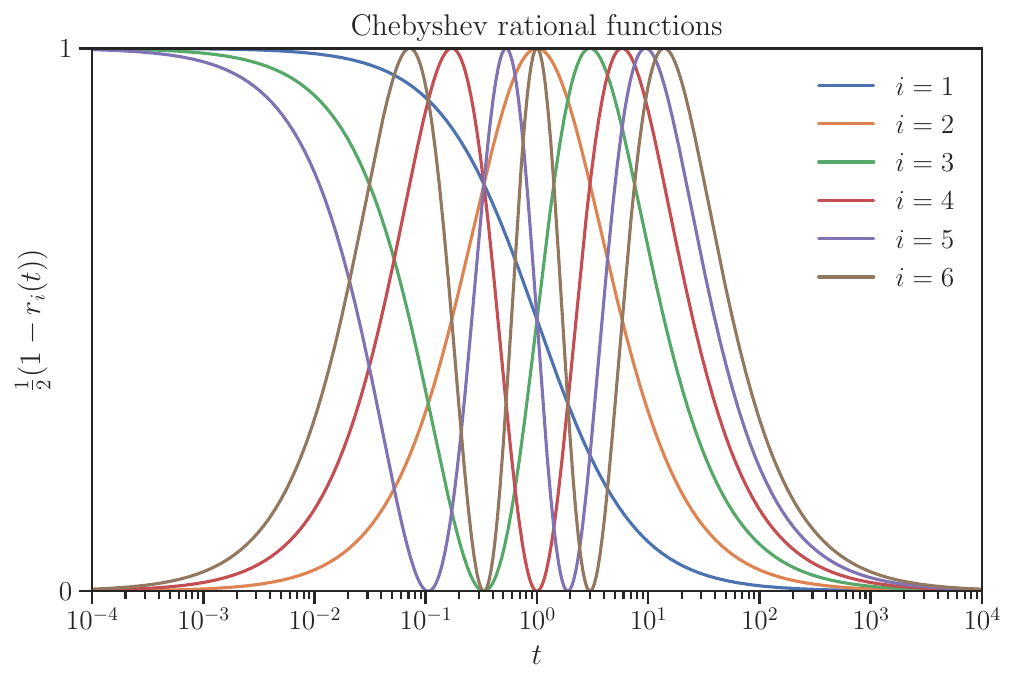}
    \caption{Chebyshev rational functions (excluding \(r_0\)) used in \eqref{eq:cheb-interp} in the logarithmic scale of \(\t\).}
    \label{fig:cheb}
\end{figure}

\Cref{fig:cheb} shows the Chebyshev rational basis functions in the form that are used on the right-hand side of \eqref{eq:cheb-interp}. These basis functions converge at \(\t \ll 1\) and \(\t \gg 1\), whereas the main variability of these functions is mostly observed near \(\mathcal{O}(t) \approx 1\). Thus, it is desirable to shift the interval of interpolation to the vicinity of \(\mathcal{O}(\t) = 1\), which can be achieved by setting the scale parameter \(\alpha\). One approach to find an optimal interpolation parameter, \(\alpha\), is to minimize the curvature of the interpolating function, which is a common practice, for instance, in smoothing splines \citep{NEWBERY-1991}. To this end, let \(w_i(\alpha)\) denote the solved weights for a given \(\alpha\). For simplicity, we transform the graph \((t, \tilde{\tau}_p)\) to \((x, y_{\alpha})\) where \(x \coloneqq (t - \alpha) / (t + \alpha)\) and \(y_{\alpha}(x) \coloneqq \sum_{i=1}^{q+1} \frac{w_i(\alpha)}{2} (1-T_i(x))\). Then, an optimal \(\alpha^*\) can be sought to minimize the arc integral of curvature squared of \(y_{\alpha}(x)\) by
\begin{equation}
    \alpha^* = \argmin_{\alpha} \int_{-1}^1 \frac{\vert y_{\alpha}''(x) \vert^2}{(1 + \vert y_{\alpha}'(x) \vert^2)^{\frac{5}{2}}} \mathrm{d} x.
\end{equation}
We note that in the absence of enough interpolant points, minimizing the curvature of the interpolating curve does not necessarily reduce interpolation error. However, when an adequate number of interpolant points are employed, the above approach can practically lead to a scale parameter \(\alpha\) that enhances the interpolation.

Finally, we note that unlike the interpolation scheme of \Cref{sec:interp-1} with the inverse monomial basis \eqref{eq:basis}, both the Pad\'{e} approximation of \eqref{eq:rational} and Chebyshev rational interpolation in \eqref{eq:cheb-interp} can interpolate \(\tau_p\) at negative values of \(\t\), namely in the domain \(\t > \t_{\inf}\) when \(t_{\inf} < 0\) (see \eqref{eq:tau-ineq-neg-1} and \eqref{eq:tau-ineq-neg-2}).


\section{Numerical Examples} \label{sec:example}

In \Cref{sec:soft}, we briefly introduce a software package we developed for the presented numerical algorithm. This package was used to produce the results in \Cref{sec:example-1} and \Cref{sec:interp-2}. Indeed, the source code to reproduce the results and plots in the following sections can be found on the documentation of the software package\footnote{See \url{https://ameli.github.io/imate}.}. \Cref{sec:example-1} considers the problem of marginal likelihood estimation, which considers a full rank correlation matrix, and for this we use the interpolation functions of \Cref{sec:interp-1}. \Cref{sec:example-2} considers the problem of ridge regression, which considers a singular matrix, and for this we use the rational polynomial interpolation method of \Cref{sec:interp-2}. We note that the interpolation with Chebyshev rational functions provide similar results to the orthogonalized inverse-monomials in \eqref{eq:approx} and we omit in our numerical examples for brevity.


\subsection{Software Package} \label{sec:soft}

The methods developed in this manuscript have been implemented into the python package \texttt{imate}, an implicit matrix trace estimator \citep{AMELI-2022-c}. This library estimates the determinant and trace of various functions of implicit matrices using either direct or stochastic estimation techniques and can process both dense matrices and large-scale sparse matrices. The main library of this package is written in \CC\ and NVIDIA\textsuperscript{\tiny{\textregistered}} CUDA and accelerated on both parallel CPU processors and CUDA-capable multi-GPU devices. The \texttt{imate} library is employed in the python package \texttt{glearn}, a machine learning library using Gaussian process regression \citep{AMELI-2022-b}.

\begin{lstlisting}[
    style=mystyle,
    language=Python,
    caption={A minimalistic usage of \texttt{imate.InterpolateSchatten} class},
    label={list:imate},
    float=tp]
# Install imate with ;#"\texttt{pip install imate}"#;
from numpy import logspace
import imate

# Generate a sample correlation matrix using the kernel ;#$e^{-r/0.1}$#;.
A = imate.sample_matrices.correlation_matrix(
        50, dimension=2, kernel='exponential', scale=0.1) ;# \label{ln:A} #;

# Create an interpolating object for
# ;#$f_p : t \mapsto \| \mathbf{A} + t \mathbf{I} \|_p$, $p=-1$#;.
f = imate.InterpolateSchatten(A, B=None, p=-1, ti=logspace(-4, 3, 8), kind='IMBF',
                                 options={'method': 'cholesky'}) ;# \label{ln:tau} #;

# Interpolate 1000 points in ;#$[10^{-4}, 10^{3}]$#;.
t = logspace(-4, 3, 1000)
y = f(t) ;# \label{ln:f} #;

# Plot the interpolated normalized curve ;#$ \tau_p(t) = \| \mathbf{A} + t \mathbf{I} \|_p / \| \mathbf{I} \|_p  $#;, compare with exact values.
f.plot(t, compare=True, normalize=True)  ;# \label{ln:plot} #;
\end{lstlisting}

\begin{figure}[bt!]
    \centering
    \includegraphics[width=0.85\textwidth]{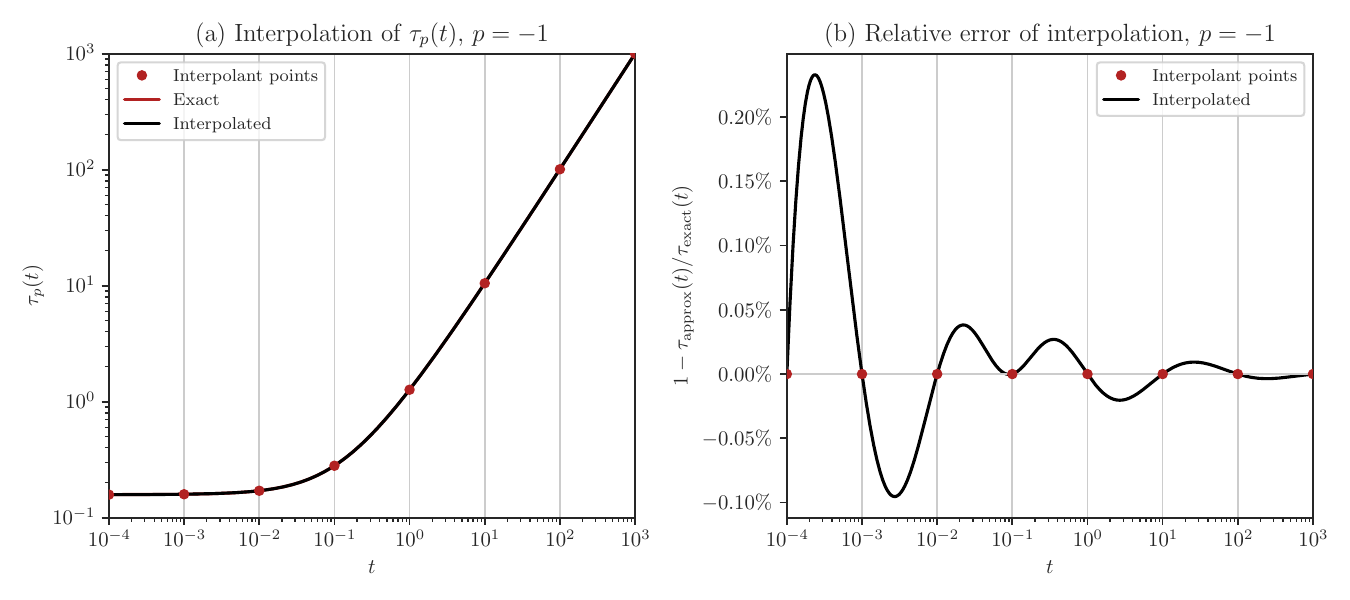}
    \caption{Result of the code in \Cref{list:imate}. (a) Comparison of interpolated versus exact value of the function \(\tau_{-1}(\t)\). The exact function (red curve) is overlaid by the interpolated curve (black curve). (b) Relative error of the comparison.}
    \label{fig:imate_mwe}
\end{figure}

In \Cref{list:imate}, we demonstrate a minimalistic usage of \texttt{imate.InterpolateSchatten} class that interpolates \(f_p: t \mapsto \| \tens{A} + t \tens{B} \|_p\). Briefly, \Cref{ln:A} generates a sample correlation matrix \(\tens{A} \in \mathcal{M}_{n, n}(\mathbb{R})\) on a randomly generated set of \(n=50^2\) points using an exponential decay kernel. In \Cref{ln:tau}, we create an instance of the class \texttt{imate.InterpolateSchatten}. Setting \texttt{B=None} indicates \(\tens{B}\) is the identity matrix using an efficient implementation that does not require storing identity matrix. The instantiation in \Cref{ln:f} internally computes \(\tau_{-1, i} = \tau_{-1}(\t_i)\) on eight interpolant points \(\t_i = 10^{-4}, 10^{-3}, \dots, 10^{3}\) and obtains the interpolation coefficients for the orthogonalized inverse-monomial basis functions \eqref{eq:basis} since \texttt{kind='IMBF'} was specified. Other possible methods can be the exact method with no interpolation (\texttt{EXT}), eigenvalue method (\texttt{EIG}) given in \eqref{eq:tr-AI}, monomial basis functions (\texttt{MBF}) given in \eqref{eq:approx-2}, Pad\'{e} rational polynomial functions (\texttt{RPF}) given in \eqref{eq:rational}, Chebyshev rational functions (\texttt{CRF}) given in \eqref{eq:cheb-interp}, or radial basis functions (\texttt{RBF}) (which we do not cover herein for brevity). The evaluation of \(\tau_{-1, i}\) can be configured by passing a dictionary of settings to the \texttt{options} argument, and we refer the interested reader to the package documentation for further details. In this minimalistic example, we compute \(\tau_{-1, i}\) using Cholesky decomposition, as further detailed in \Cref{sec:example-1}. Other methods include stochastic Lanczos quadrature or Hutchinson estimation; we compare such methods in \Cref{sec:example-2}. Once the interpolation object is initialized, future calls to interpolate an arbitrary number of points \(\t\) are returned almost instantly. In \Cref{ln:f}, the interpolation is performed on \(1000\) points in the interval \([10^{-4}, 10^3]\) spaced uniformly on the logarithmic scale. A comparison of the interpolated result versus the exact solution are shown in \Cref{fig:imate_mwe}. It can be seen that with only eight interpolant points, the relative error of interpolation over a wide range of parameter \(\t\) is around \(0.1\%\). 

\subsection{Marginal Likelihood Estimation for Gaussian Process Regression} \label{sec:example-1}

Here we generate a full rank correlation matrix from a spatially correlated set of points \(\vect{x} \in \mathcal{D} = [0,1]^2\). To define a spatial correlation function, we use the isotropic exponential decay kernel given by
\begin{equation}
    \kappa(\vect{x},\vect{x}'\vert \rho) = \exp \left(-\frac{\| \vect{x} - \vect{x}' \|_2}{\rho} \right), \label{eq:exp-decay}
\end{equation}
where \(\rho\) is the correlation scale, set to \(\rho = 0.1\). The above exponential decay kernel represents an Ornstein-Uhlenbeck random process, which is a Gaussian and zeroth-order Markov process \cite[p. 85]{RASMUSSEN-2006}. To produce discrete data, we sample \(n = 50^2\) points from \(\mathcal{D}\), which yields the symmetric and positive-definite correlation matrix \(\tens{A}\) with the components \(A_{ij} = \kappa(\vect{x}_i,\vect{x}_j\vert \rho)\). We aim to interpolate functions
\begin{subequations}
\begin{align}
    &\log \det\left( \tens{A} + \t \tens{I} \right) = n \log \tau_0(\t) , \label{eq:tau-K-1} \\
    &\trace\left( (\tens{A} + \t \tens{I})^p \right) = n(\tau_p(\t))^{p}, \label{eq:tau-K-2}
\end{align}
\end{subequations}
for \(p = -1, -2\), which appear in many statistical applications, such as the estimation of noise in Gaussian process regression~\citep{AMELI-2022-a}. Specifically, the above functions for \(p = 0\), \(-1\), and \(-2\) appear in the corresponding likelihood function, and its Jacobian and Hessian, respectively.

We compute the exact value of \(\tau_p(\t)\) for \(p \in \mathbb{Z}_{\leq 0}\) (either at interpolant points \(\t_i\) or at all points \(\t\) for the purpose of benchmark comparison) as follows. We compute the Cholesky factorization of \( (\tens{A} + t \tens{I})^{\vert p \vert} = \tens{L}_{\vert p \vert} \tens{L}_{\vert p \vert}^{\intercal}\), where \(\tens{L}_{\vert p \vert}\) is lower triangular. Then
\begin{subequations}
\begin{align}
    &\log \det(\tens{A} + t \tens{I}) =  2 \sum_{i=1}^n \log ((\tens{L}_{1})_{ii}), \label{eq:det-chol} \\
    &\trace \left( (\tens{A} + t \tens{I})^p \right) = \trace(\tens{L}_{\vert p \vert}^{-\intercal} \tens{L}_{\vert p \vert}^{-1}) 
    = \trace(\tens{L}_{\vert p \vert}^{-1} \tens{L}_{\vert p \vert}^{-\intercal})
    = \| \tens{L}_{\vert p \vert}^{-1} \|^2_{F}, \quad p \in \mathbb{Z}_{< 0},
    \label{eq:trace-chol}
\end{align}
\end{subequations}
where \((\tens{L}_{1})_{ii}\) is the \(i\)\textsuperscript{th} diagonal element of \(\tens{L}_{1}\) and \(\| \cdot \|_F\) is the Frobenius norm. The second equality in \eqref{eq:trace-chol} employs the cyclic property of the trace operator. A simple method to compute \(\| \tens{L}_{\vert p \vert}^{-1} \|^2_F\) without storing \(\tens{L}_{\vert p \vert}^{-1}\) is to solve the lower triangular system \(\tens{L}_{\vert p \vert} \vect{x}_i = \vect{e}_i\) for \(\vect{x}_i\), \(i = 1,\dots,n\), where \(\vect{e}_i = (0, \dots, 0,1,0,\dots,0)^{\intercal}\) is a column vector of zeros, except, its \(i\)\textsuperscript{th} entry is one. The solution vector \(\vect{x}_i\) is the \(i\)\textsuperscript{th} column of \(\tens{L}_{\vert p \vert}^{-1}\). Thus, \(\| \tens{L}_{\vert p \vert}^{-1} \|^2_{F} = \sum_{i = 1}^n \| \vect{x}_i \|^2\). This method is memory efficient since the vectors \(\vect{x}_i\) do not need to be stored.

We note that the complexity of the interpolation method is the number of evaluations of \(\tau_p\) at interpolant points \(t_i\) and at \(t=0\) (which is proportional to \(q\)) times the complexity of computing \(\tau_{p}\) at a single point \(t\). For instance, by using the Cholesky method in \eqref{eq:det-chol} or \eqref{eq:trace-chol} which costs \(\mathcal{O}(\frac{1}{3}n^3)\) for a matrix of size \(n\), the complexity of the interpolation method is \(\mathcal{O}(\frac{1}{3} n^3 q)\).

\begin{remark}[Case of Sparse Matrices]
There exist efficient methods to compute the Cholesky factorization of sparse matrices (see \eg \citet[Ch. 4]{DAVIS-2006}). Also, the inverse of the sparse triangular matrix \(\tens{L}_{\vert p \vert}\) can be computed at \(\mathcal{O}(n^2)\) complexity \cite[pp. 93-95]{STEWART-1998}, and a linear system with both sparse kernel \(\tens{L}_{\vert p \vert}\) and sparse right-hand side \(\vect{e}_i\) can be solved efficiently (see \citet[Sec. 3.2]{DAVIS-2006}). 
\end{remark}

\begin{figure}[hpt]
    \centering
    \subfloat{\includegraphics[width=\textwidth]{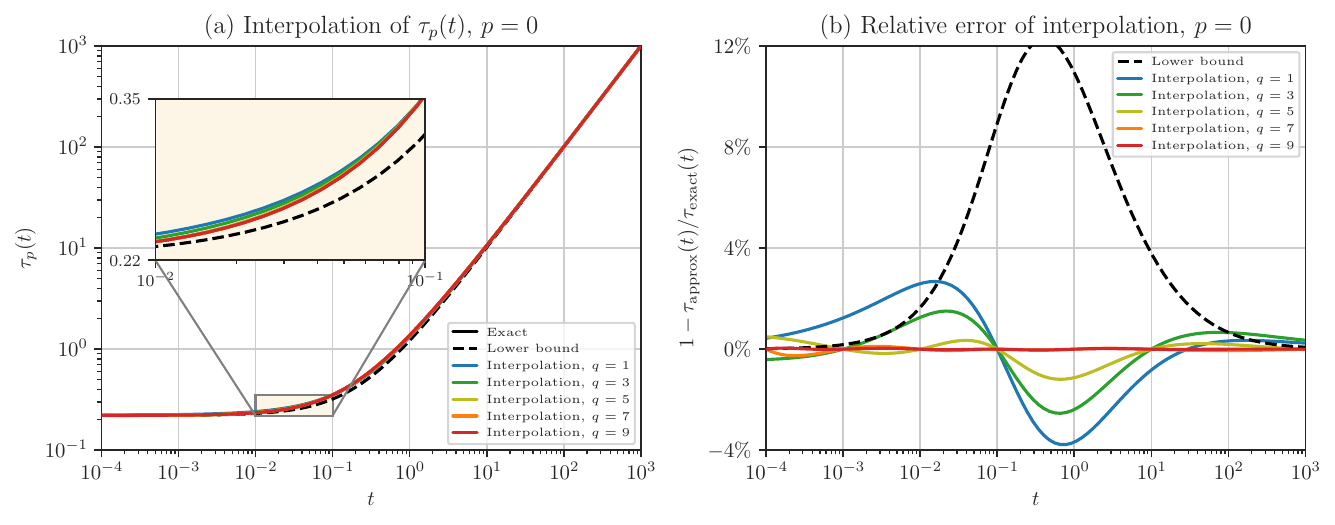}}
    \vspace{-5.5mm} \\
    \subfloat{\includegraphics[width=\textwidth]{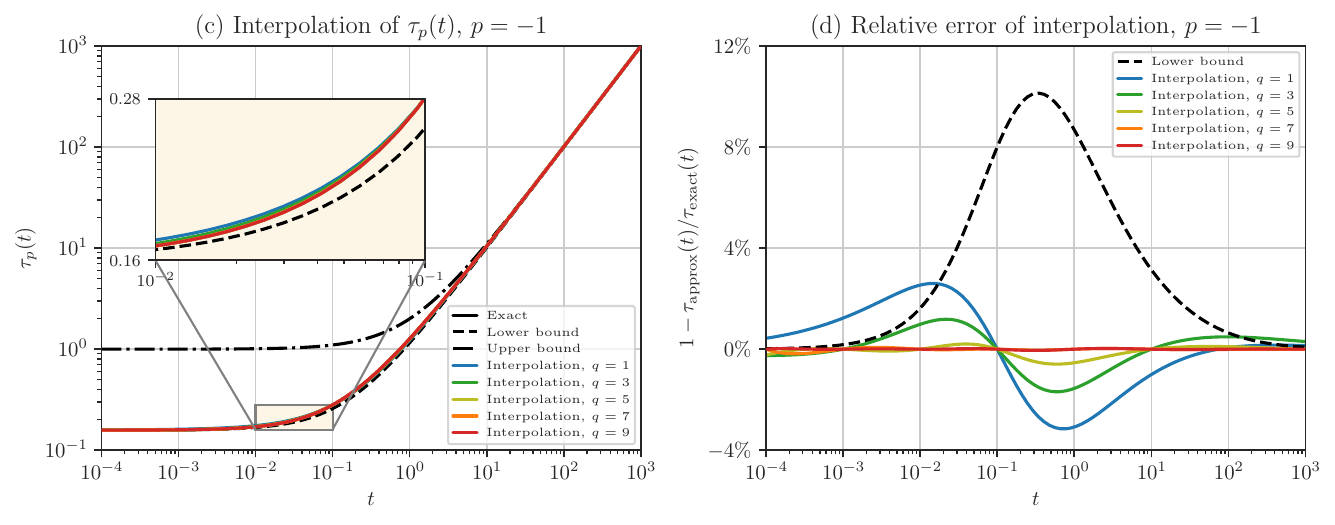}}
    \vspace{-5.5mm} \\
    \subfloat{\includegraphics[width=\textwidth]{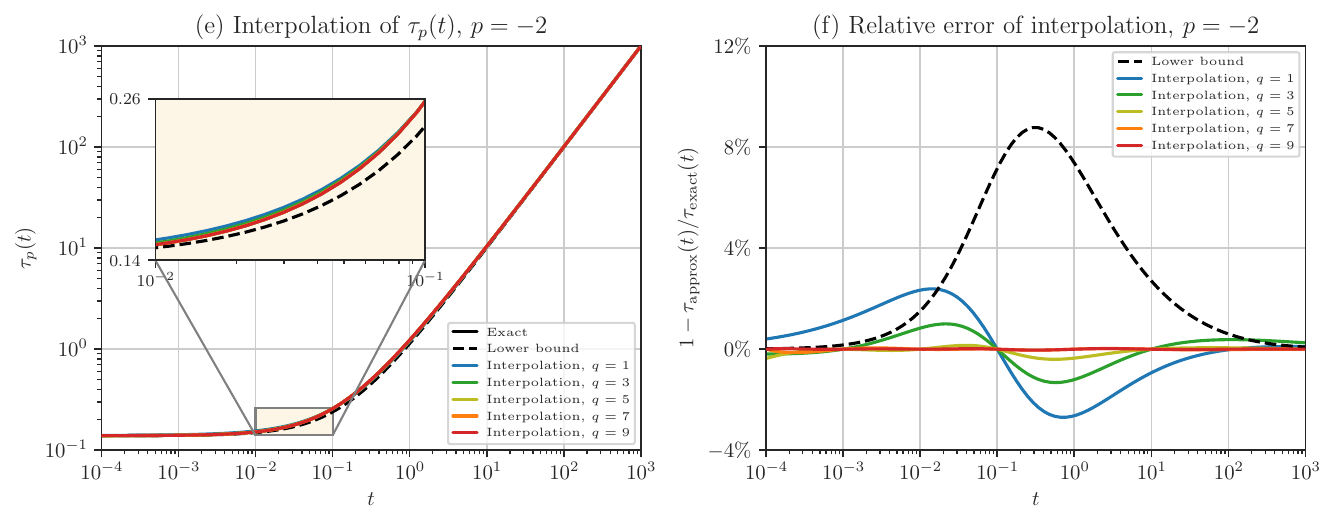}}
    \caption{Left columns: Comparison of the exact function \(\tau_p(\t)\), bounds \(\hat{\tau}_p\), and the interpolations \(\tilde{\tau}_p(\t)\) for various numbers of interpolant points. The interpolation becomes almost indistinguishable from the exact solution once 5 or more interpolation points are used. Right columns: Relative error of the interpolations and the bounds. The interpolations using 7 and 9 points lead to relative errors of less than \(0.02 \%\) and \(0.01 \%\), respectively. Rows correspond to \(p=0\), \(-1\), and \(-2\), respectively.}
    \label{fig:trace}
\end{figure}

The exact value of \(\tau_p(\t)\), for \(p = 0, -1, -2\), computed directly using the Cholesky factorization method described above are respectively shown in \Cref{fig:trace}(a, c, e) by the solid black curve (overlaid by the red curve) in the range \(\t \in [10^{-4},10^3]\). The dashed black curves in \Cref{fig:trace}(a, c, e) are the lower bounds \(\hat{\tau}_p(\t)\) given by \eqref{eq:tau-upper}, which can be thought of as the estimation with zero interpolant points, \ie \(q = 0\). For completeness, we have also shown the upper bound of \(\tau_{-1}(\t)\) by the black dash-dot line in \Cref{fig:trace}(c), given by
\begin{equation}
    \check{\tau}_{-1}(\t) \coloneqq 1+\t \geq \tau_{-1}(\t). \label{eq:lower}
\end{equation}
The above upper bound can be obtained from \eqref{eq:trace-lower-bound} and the fact that \(\trace(\tens{A}) = n\), since the diagonals of the correlation matrix are \(1\). However, unlike the lower bound in \eqref{eq:anti-ineq}, the upper bound \eqref{eq:lower} is not useful for approximation as it does not asymptote to \(\tau_{-1}(\t)\) at small \(\t\). Nonetheless, both the lower and upper bounds asymptote to \(\t\) at large \(\t\).

To estimate \(\tau_p\), we used the interpolation function in \eqref{eq:approx} with the set of orthonormal basis functions in \Cref{table:ortho-funct}.  The colored solid lines in \Cref{fig:trace}(a, c, e) are the interpolations \(\tilde{\tau}_{p}(\t)\) with \(q = 1,3,5,7\), and \(9\) interpolant points, \(\t_i\), spanning from \(10^{-4}\) to \(10^3\).  It can be seen from the embedded diagrams in \Cref{fig:trace}(a, c, e) that \(\tilde{\tau}_p(\t)\) is remarkably close to the true function value. In practice, fewer interpolant points in a small range, \eg \([10^{-2},10^2]\), are sufficient to effectively interpolate \(\tau_p\).

To better compare the exact and interpolated functions, the relative error of the interpolations is shown in \Cref{fig:trace}(b, d, f). The relative error of the lower bound (dashed curve) rapidly vanishes at both ends, namely, at \(\t \ll \tau_{p, 0}\) and \(\t \gg \tau_{p, 0}\), where \(\tau_{0, 0} = 0.22\), \(\tau_{-1, 0} = 0.16\), and \(\tau_{-2, 0} = 0.14\). The absolute error of the upper bound is highest at \(\mathcal{O}(\t \tau_{p, 0}^{-1}) = 1\), or \(\t \approx \tau_{p, 0}\), which is slightly to the left of the relative error peak on each diagram.

Based on the lower bound, we distribute the interpolant points, \(\t_i\), almost evenly around \(\t \approx \tau_{p, 0}\) where the lower bound has the highest error. The blue curve in \Cref{fig:trace}(b, d, f) corresponds to the case with only one interpolation point at \(\t_1 = 10^{-1}\), which already leads to a relative error less than \(3 \%\) almost everywhere. On the other hand, with only nine interpolation points \(\t_i \in \{10^{-4},4\times 10^{-4},10^{-3},10^{-2},\dots,10^3\}\) the relative error becomes less than \(0.01 \% \). Beyond the strong accuracy shown by the relative errors, the absolute errors are more compelling since \( \tau_p(\t)\) decays by orders of magnitude at large \(\t\), making the absolute error negligible at \(\t \gg \tau_{p, 0}\).


\subsection{Ridge Regression with Generalized Cross-Validation} \label{sec:example-2}

Here we calculate the optimal regularization parameter for a linear ridge regression model using generalized cross-validation (GCV). Consider the linear model \(\vect{y} = \tens{X} \vect{\beta} + \vect{\epsilon}\), where \(\vect{y} \in \mathbb{R}^n\) is a column vector of given data, \(\tens{X} \in \mathcal{M}_{n, m}(\mathbb{R})\) is the known design matrix representing \(m\) basis functions where \(m < n\), \(\vect{\beta} \in \mathbb{R}^{m}\) is the unknown coefficients of the linear model, and \(\vect{\epsilon} \sim \mathcal{N}(\vect{0},\sigma^2 \tens{K})\) is the correlated residual error of the model, which is a zero-mean Gaussian random vector with the symmetric and positive-definite correlation matrix \(\tens{K}\) and unknown variance \(\sigma^2\). A generalized least-squares solution to this problem minimizes the square Mahalanobis distance \(\| \vect{y} - \tens{X} \vect{\beta} \|_{\tens{K}^{-1}}^2 \coloneqq (\vect{y} - \tens{X} \vect{\beta})^{\intercal} \tens{K}^{-1} (\vect{y} - \tens{X} \vect{\beta})\) yielding an estimation of \(\vect{\beta}\) by \(\hat{\vect{\beta}} = (\tens{X}^{\intercal} \tens{K}^{-1} \tens{X})^{-1} \tens{X}^{\intercal} \tens{K}^{-1} \vect{y}\) \citep[p. 67]{SEBER-2012}.

When \(\tens{X}\) is not full rank, the least-squares problem is not well-conditioned. A resolution of the ill-conditioned problems is the ridge (Tikhonov) regularization, where the function \(\| \vect{y} - \tens{X} \vect{\beta} \|_{\tens{K}^{-1}}^2 + n\theta \| \vect{\beta} \|_{\gtens{\Omega}}^2\) is minimized instead \citep[Sec. 12.5.2]{SEBER-2012}. Here, the penalty term is \(\Vert \vect{\beta} \Vert_{\gtens{\Omega}}^2 = \vect{\beta}^{\intercal} \gtens{\Omega} \vect{\beta}\) where \(\gtens{\Omega}\) is the symmetric and positive-definite penalty matrix. The estimate of \(\vect{\beta}\) using the penalty term becomes
\begin{equation}
    \hat{\vect{\beta}} = (\tens{X}^{\intercal} \tens{K}^{-1} \tens{X} + n\theta \gtens{\Omega})^{-1} \tens{X}^{\intercal} \tens{K}^{-1} \vect{y}. \label{eq:beta}
\end{equation}
Also, the fitted values on the training points are \(\hat{\vect{y}} = \tens{X} \hat{\vect{\beta}}\), which can be written as \(\hat{\vect{y}} = \tens{S}_{\theta} \vect{y}\), where the smoother matrix \(\tens{S}_{\theta}\) is defined by
\begin{equation}
    \tens{S}_{\theta} \coloneqq \tens{X} (\tens{X}^{\intercal} \tens{K}^{-1} \tens{X} + n\theta \gtens{\Omega})^{-1} \tens{X}^{\intercal} \tens{K}^{-1}. \label{eq:smoother}
\end{equation}

The regularization parameter, \(\theta\), plays a crucial role to balance the residual error versus the added penalty term. The generalized cross-validation method  \citep{WAHBA-1977,CRAVEN-1978,GOLUB-1979} is a popular way to seek an optimal regularization parameter without needing to estimate the error variance \(\sigma^2\). Namely, the regularization parameter is sought as the minimizer of
\begin{equation}
    V(\theta) \coloneqq \frac{\frac{1}{n} \left\| (\tens{I} - \tens{S}_{\theta} ) \vect{y} \right\|_{\tens{K}^{-1}}^2}{\left( \frac{1}{n} \trace( \tens{I} - \tens{S}_{\theta}) \right)^2}, \label{eq:gcv}
\end{equation}
\citep[p. 244]{HASTIE-2001}\footnote{
The function \eqref{eq:gcv} is modified to incorporate the correlation of error using \(\tens{K}\), and can be derived from the conventional definition of generalized cross-validation function for the decorrelated error \(\vect{\epsilon}' \coloneqq \tens{L}^{-1} \epsilon \sim \mathcal{N}(\vect{0}, \sigma^2 \tens{I})\) where \(\tens{K} = \tens{L} \tens{L}^{\intercal}\) is the Cholesky decomposition of \(\tens{K}\).}. For large matrices, it is difficult is to compute \(\trace(\tens{S}_{\theta})\) (also known as the effective degrees of freedom) in the denominator of \eqref{eq:gcv}, and several methods have been developed to address this problem \citep{GOLUB-1997}, \citep{LUKAS-2010}. 

\subsubsection{Estimating the Trace}
Using the presented interpolation method, we aim to estimate \(\trace(\tens{S}_{\theta})\). Let \(\gtens{\Omega} = \tens{L} \tens{L}^{\intercal}\) be the Cholesky decomposition of \(\gtens{\Omega}\). Using the cyclic property of trace operator, we have
\begin{align}
    \trace(\tens{S}_{\theta}) &=\trace ( (\tens{X}^{\intercal} \tens{K}^{-1} \tens{X} + n\theta \gtens{\Omega})^{-1} \tens{X}^{\intercal} \tens{K}^{-1} \tens{X} ) \notag \\
                              &= \trace ( \tens{I}_{m \times m} - n\theta (\tens{X}^{\intercal} \tens{K}^{-1} \tens{X} + n\theta \gtens{\Omega})^{-1} \gtens{\Omega} ) \notag \\
                              &= m - n\theta \trace ( \tens{L}^{\intercal} ( \tens{X}^{\intercal} \tens{K}^{-1} \tens{X} + n\theta \gtens{\Omega})^{-1} \tens{L} ) \notag \\
                              &= m - n\theta \trace ( (\tens{L}^{-1} \tens{X}^{\intercal} \tens{K}^{-1} \tens{X} \tens{L}^{-\intercal} + n\theta \tens{I})^{-1}). \label{eq:trick}
\end{align}
In the above, \(\tens{I}_{m \times m}\) is identity matrix of size \(m\). To compute the above term, we interpolate
\begin{equation}
    \trace\left( (\tens{A} + \t \tens{I})^{-1} \right) = m(\tau_{-1}(\t))^{-1}, \label{eq:tau-shifted}
\end{equation}
where \(\t \coloneqq n\theta - s\) and
\begin{equation*}
    \tens{A} \coloneqq \tens{L}^{-1} \tens{X}^{\intercal} \tens{K}^{-1} \tens{X} \tens{L}^{-\intercal} + s \tens{I}.
\end{equation*}
We note that the size of \(\tens{A}\) and \(\tens{I}\) is \(m\). Also, \(\tens{A}\) is symmetric and positive-definite since it can be written as a Gramian matrix. The purpose of the fixed parameter \(s \ll 1\) is to slightly shift the singular matrix \(\tens{L}^{-1} \tens{X}^{\intercal} \tens{K}^{-1} \tens{X} \tens{L}^{-\intercal}\) to make \(\tens{A}\) non-singular. The shift is necessary since without it, \eqref{eq:tau-shifted} is undefined at \(\t = 0\), and we cannot compute \(\tau_{-1, 0} = m / \trace(\tens{A}^{-1})\). Also, the shift can improve interpolation by relocating the origin of \(\t\) to the vicinity of the interval where we are interested to compute \(V(\theta)\).

For simplicity in our numerical experiment, we set \(\tens{K}\) and \(\gtens{\Omega}\) to identity matrices of sizes \(n\) and \(m\), respectively. We also set \(s = 10^{-3}\). We create an ill-conditioned design matrix \(\tens{X}\) for our numerical example using singular value decomposition \(\tens{X} = \tens{U} \gtens{\Sigma} \tens{V}^{\intercal}\). The orthogonal matrices \(\tens{U} \in \mathcal{M}_{n, n}(\mathbb{R})\) and \(\tens{V} \in \mathcal{M}_{m, m}(\mathbb{R})\) were produced by the Householder matrices
\begin{equation*}
    \tens{U} \coloneqq \tens{I} - 2 \frac{\vect{u} \vect{u}^{\intercal}}{\| \vect{u} \|_2^2},
    \qquad \text{and} \qquad 
    \tens{V} \coloneqq \tens{I} - 2 \frac{\vect{v} \vect{v}^{\intercal}}{\| \vect{v} \|_2^2},
\end{equation*}
where \(\vect{u}\in \mathbb{R}^n\) and \(\vect{v} \in \mathbb{R}^m\) are random vectors (see also \citep[Sec. 10]{GOLUB-1997}). The diagonal matrix \(\gtens{\Sigma} \in \mathcal{M}_{n, m}(\mathbb{R})\) was defined by 
\begin{equation}
    \Sigma_{ii} \coloneqq \exp \bigg( -40 \Big( \frac{i-1}{m}\Big)^{\sfrac{3}{4}} \bigg), \quad i = 1,\dots,m. \label{eq:singulars}
\end{equation}
We set \(n = 10^3\) and \(m = 500\). We generated data by letting \(\vect{y} = \tens{X} \vect{\beta} + \vect{\epsilon}\), where \(\vect{\beta}\), and \(\vect{\epsilon}\), were randomly generated with a unit variance, and \(\sigma = 0.4\), respectively.

\begin{figure}[t!]
    \centering
    \includegraphics[width=\textwidth]{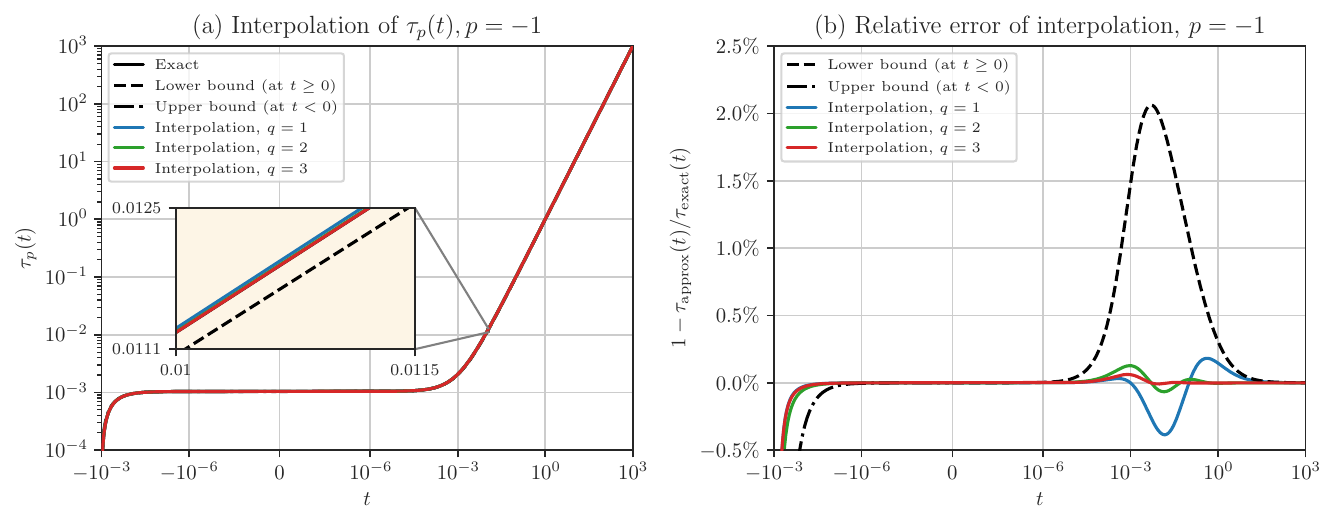}
    \caption{
    (a) The exact solution \(\tau_{-1}(\t)\), bounds \(\hat{\tau}_{-1}(\t)\), and Pad\'{e} rational polynomial interpolations \(\tilde{\tau}_{-1}(\t)\) for \(q = 1, 2, 3\) are shown. The green curve and the exact solution in the solid black curve are overlaid by the red curve. The embedded diagram (with linear axes) magnifies a portion of the curves with the highest interpolation error. (b) The relative error of the curves in (a) with respect to the exact solution is shown. In both diagrams (a) and (b), the horizontal axis in the interval \([-10^{-6},10^{-6}]\) is linear, but outside this interval, the axis is logarithmic.}
    \label{fig:trace-2}
\end{figure}

We computed the exact solution of \(\tau_{-1}(\t)\) in \eqref{eq:tau-shifted}, and interpolation points \(\tau_{-1}(\t_i)\), using the Cholesky factorization method described by \eqref{eq:trace-chol}. The exact solution is shown by the solid black curve in \Cref{fig:trace-2} (overlaid by the red curve) with \(\tau_{-1, 0} = 960.5^{-1}\). The lower bound \(\hat{\tau}_{-1}(\t)\) from \eqref{eq:anti-ineq} is shown by the dashed black curve for \(\t > 0\). In contrast, at \(\t \in (\t_{\inf},0]\), the upper bound from \eqref{eq:anti-ineq-rev} is shown, where \(\t_{\inf} = -\min(\eig(\tens{A}))\)  and \(\min(\eig(\tens{A})) \approx s = 10^{-3}\) is the smallest eigenvalue of \(\tens{A}\). The relative error of the bounds with respect to the exact solution are shown in \Cref{fig:trace-2}(b). The peak of the absolute error of the lower bound is located approximately at \(\t \approx \tau_{-1, 0} \approx 10^{-3}\), and the peak of its relative error of the lower bound is slightly to the right of this value.

We sought to interpolate \eqref{eq:tau-shifted} in the interval \(\theta \in [10^{-7},10]\). Accordingly, since we set \(s = 10^{-3}\), we shifted the origin of \(\t = n \theta - s\) inside the interval \(n \theta \in [10^{-4},10^4]\). Thus, we approximately had \(-10^{-3} < \t < 10^4\). Because this interval contains the origin, we employed the Pad\'{e} rational polynomial interpolation method in \Cref{sec:interp-2}. (Recall that at small \(\t\), particularly at \(\vert \t \vert \ll \tau_{-1, 0}\), the rational polynomial interpolation performs better than the basis functions interpolation.) We distributed the interpolant points at \(\t_i \geq \tau_{-1, 0} \approx 10^{-3}\) where the rational polynomial interpolation has to adhere to the exact solution. 

The interpolation function \(\tilde{\tau}_{-1}(\t)\) with \(q = 1, 2, 3\) is shown in \Cref{fig:trace-2}(a) using \(2q\) interpolation points \(t_i\) in the interval \(t_i \in [5 \times 10^{-3}, 5]\) that are equally distanced in the logarithmic scale. The red curve corresponding to \(q=3\) and the black curve (exact solution) are indistinguishable even in the embedded diagram that magnifies the location with the highest error. The relative error of the interpolations is shown by \Cref{fig:trace-2}(b). On the far left of the range of \(\t\), the error spikes due to the singularity of the matrix \(\tens{X}\), which makes \(\tau_{-1}(\t)\) undefined at \(\t = - 10^{-3}\), corresponding to \(\theta = 0\). On the rest of the range, the green and red curves respectively show less than \(0.1 \%\) and \( 0.05 \%\) relative errors, which are compelling accuracy for a broad range of \(\t\), and achieved with only four and six interpolation points, respectively.

\subsubsection{Optimization of Generalized Cross-Validation}

\begin{figure}[bt!]
    \centering
    \includegraphics[width=0.6\textwidth]{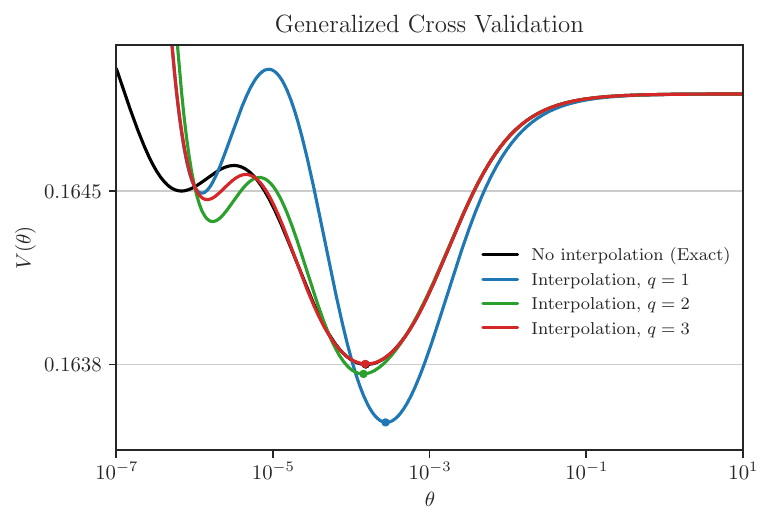}
    \caption{The generalized cross-validation function is shown, where the black and colored curves correspond respectively to the exact and interpolated computation of \(\tau_{-1}(\t)\) in the denominator of \(V(\theta)\). The global minimum of each curve  is shown by a dot.}
    \label{fig:gcv}
\end{figure}

Here we apply the result of our trace interpolations above to solve the generalized cross-validation problem. The function \(V(\theta)\) from \eqref{eq:gcv} is plotted in \Cref{fig:gcv}, with the black curve, corresponding to the exact solution with \(\tau_{-1}(\t)\) applied in the denominator of \(V(\theta)\), serving as a benchmark for comparison. The blue, green, and red curves correspond to the proceeding trace interpolations applied in the denominator of \(V(\theta)\). The interpolated curves exhibit both local minima of \(V(\theta)\) similar to the benchmark curve, but with slight differences in the positions of the minima. Due to the singularity at \(\theta = 0\), the interpolations of \(\tau_{-1}(\t)\) become less accurate at low values of \(\theta\). At higher values of \(\theta\), all curves steadily asymptote to a constant. We note that the results in \Cref{fig:gcv} are compelling since the estimation of \(V(\theta)\) is sensitive to the interpolation of its denominator. Namely, a consistent interpolation accuracy over all the parameter range is essential to capture the qualitative shape of \(V(\theta)\) correctly.

\begin{table}[t]
    \centering
    \caption{
    Comparison of methods to optimize the regularization parameter \(\theta\), with and without interpolation of \(\tau_{-1}(\t)\), and by various algorithms of computing trace of a matrix inverse.
    } \label{table:gcv}
    \begin{adjustbox}{width=1\textwidth}
    \begin{small}
    \begin{tabular}{l l c c r r c c c r r}
        \toprule
        \multicolumn{2}{c}{Computing \(\tau_{-1}(\t)\)} & \multicolumn{2}{c}{Iterations} & \multicolumn{2}{c}{Time (sec)} & \multicolumn{2}{c}{Results} & \multicolumn{3}{c}{Relative Error} \\
        \cmidrule(lr){1-2} \cmidrule(lr){3-4} \cmidrule(lr){5-6} \cmidrule(lr){7-8} \cmidrule(lr){9-11}
        Algorithm & Interpolation Method & \(N_{\text{tr}}\) & \(N_{\text{tot}}\) & \multicolumn{1}{c}{\(T_{\text{tr}}\)} & \multicolumn{1}{c}{\(T_{\text{tot}}\)} & \(V(\theta^*)\) & \(\log_{10} \theta^{\ast}\) & \(\frac{\vert \Delta \log_{10} \theta^{\ast}\vert}{\vert \log_{10} \theta^{\ast} \vert}\) & \multicolumn{1}{c}{\(\frac{\Vert \Delta \hat{\vect{\beta}} \Vert_2}{\Vert \hat{\vect{\beta}}\Vert_2}\)} & \multicolumn{1}{c}{\(\frac{\Vert \Delta \hat{\vect{y}} \Vert_2}{\Vert \hat{\vect{y}} \Vert_2}\)} \\
        \midrule
        Cholesky & No interpolation (exact)  & 282 & 282 &  27.49 & 30.90 & 0.16376 & -3.8164 & 0.00 \% &  0.00 \% &  0.00 \% \\
        & Rational polynomial, \(q = 1\)     & 3   & 364 &   0.29 &  4.70 & 0.16352 & -3.5628 & 6.65 \% & 29.71 \% & 17.59 \% \\
        & Rational polynomial, \(q = 2\)     & 5   & 282 &   0.49 &  3.93 & 0.16372 & -3.8446 & 0.74 \% &  3.69 \% &  1.95 \% \\
        & Rational polynomial, \(q = 3\)     & 7   & 284 &   0.69 &  4.12 & 0.16376 & -3.8218 & 0.14 \% &  0.71 \% &  0.37 \% \\
        \rule{0pt}{0ex} \\

        Hutchinson & No interpolation        & 312 & 312 &  61.33 & 65.14 & 0.16372 & -3.7939 & 0.59 \% &  2.98 \% &  1.58 \% \\
        & Rational polynomial, \(q = 1\)     & 3   & 364 &   0.57 &  4.96 & 0.16348 & -3.5649 & 6.59 \% & 29.49 \% & 17.45 \% \\
        & Rational polynomial, \(q = 2\)     & 5   & 282 &   0.88 &  4.29 & 0.16371 & -3.8274 & 0.29 \% &  1.45 \% &  0.76 \% \\
        & Rational polynomial, \(q = 3\)     & 7   & 282 &   1.25 &  4.66 & 0.16374 & -3.8119 & 0.12 \% &  0.61 \% &  0.32 \% \\
        \rule{0pt}{0ex} \\

        SLQ & No interpolation               & 322 & 322 &  66.67 & 88.75 & 0.16373 & -3.7939 & 0.59 \% &  2.98 \% &  1.58 \% \\
        & Rational polynomial, \(q = 1\)     & 3   & 364 &   0.58 &  5.04 & 0.16352 & -3.5597 & 6.73 \% & 30.03 \% & 17.81 \% \\
        & Rational polynomial, \(q = 2\)     & 5   & 282 &   1.03 &  4.52 & 0.16376 & -3.8778 & 1.61 \% &  7.88 \% &  4.18 \% \\
        & Rational polynomial, \(q = 3\)     & 7   & 282 &   0.70 &  4.13 & 0.16378 & -3.7770 & 1.03 \% &  5.17 \% &  2.76 \% \\
    \bottomrule
    \end{tabular}
    \end{small}
    \end{adjustbox}
\end{table}

The global minimum of \(V(\theta)\) at \(\theta = \theta^*\) is the optimal compromise between an ill-conditioned regression problem (small \(\theta)\) and a highly regularized regression problem (large \(\theta\)). We aimed to test the practicality of our interpolation method in searching the global minimum, \(V(\theta^*)\), by numerical optimization. We note that we constructed \(\tens{X}\) in \eqref{eq:singulars} so that \(V(\theta)\) would have two local minima thus making optimization less trivial. In general, the generalized cross-validation function may have more than one local minimum necessitating global search algorithms \citep{KENT-2000}. The optimization was performed using a differential evolution optimization method \citep{STORN-1997} with a best/1/exp strategy and 40 initial guess points. The results are shown in the first four rows of \Cref{table:gcv}, where the trace of a matrix inverse is computed by the Cholesky factorization described in \eqref{eq:trace-chol}. In the first row, \(\tau_{-1}\) is computed exactly, \ie without interpolation, at all requested locations \(\t\) during the optimization procedure. On the second to fourth rows, \(\tau_{-1}\) is first pre-computed at the interpolation points, \(\t_i\), by the Cholesky factorization, and then the interpolation is subsequently used during the optimization procedure.

In the table, \(N_{\text{tr}}\) counts the number of exact evaluations of \(\tau_{-1}\). For the Pad\'{e} rational polynomial interpolation method of degree \(q\), we had \(N_{\text{tr}} = 2q + 1\), accounting for \(2q\) interpolant points in addition to the evaluation of \(\tau_{-1, 0}\) at \(\t = 0\). \(N_{\text{tot}}\) is the total number of estimations of \(\tau_{-1}\) during the optimization process. In the first row, \(N_{\text{tr}} = N_{\text{tot}}\) as all points are evaluated exactly, \ie without interpolation. However, for the interpolation methods, \(N_{\text{tot}}\) consists of \(N_{\text{tr}}\) plus the evaluations of \(\tau_{-1}\) via interpolation.

The exact computations of \(\tau_{-1}\) (at \(N_{\text{tr}}\) points) are the most computationally expensive part of the overall process. Our numerical experiments were performed on the Intel Xeon E5-2640 v4 processor using shared memory parallelism. We measured computational costs by the total CPU processing time of all computing cores. \(T_{\text{tr}}\) denotes the processing time of computing \(\tau_{-1}\) exactly at the \(N_{\text{tr}}\) points. \(T_{\text{tot}}\) measures the processing time of the overall optimization, which includes \(T_{\text{tr}}\). As shown, the interpolation methods took significantly less processing time compared to no interpolation, namely, by two orders of magnitude for \(T_{\text{tr}}\), and an order of magnitude for \(T_{\text{tot}}\). We also observe that without interpolation, \(T_{\mathrm{tr}}\) is the dominant part of the total processing time, \(T_{\mathrm{tot}}\). In contrast, with interpolation, \(T_{\text{tr}}\) becomes so small that \(T_{\mathrm{tot}}\) is dominated by the cost of evaluating the numerator of \(V(\theta)\) in \eqref{eq:gcv}, which is proportional to \(N_{\mathrm{tot}}\).

The results of computing the optimized parameter, \(\theta^*\), and the corresponding minimum, \(V(\theta^*)\), are shown in the seventh and eighth columns of \Cref{table:gcv}, respectively. The ninth column is the relative error of estimating \(\theta^*\), and obtained by comparing \(\log_{10} \theta^{*}\) between the interpolated and benchmark solution (\ie first row). The last two columns are the \(\ell^2\) norm of the error of \(\hat{\vect{\beta}}\) (using \eqref{eq:beta}) and \(\hat{\vect{y}}\) compared to their exact solution, and their relative error are obtained by normalizing with the \(\ell^2\) norm of their exact solution. We observed, for example for \(q=2\), that with one-tenth of the processing time, an accuracy of \(2 \%\) error for \(\hat{\vect{y}}\) is achieved, which is generally sufficient in practical applications. Also, for \(q=3\), the error reduces to \(<\)1\% with similar processing time. In general, the error can be improved simply by using more interpolant points. We have found that simple heuristics for setting defaults for \(q\) are broadly effective. Namely, if \(\theta^{\ast}\) is expected to be found in an known interval, one can use a small number of interpolating points (\(q = 1 \sim 2\)) on the boundary or center of the interval. If there is no prior knowledge of the range, one can let an optimization scheme search for \(\theta^{\ast}\) in a wide logarithmic range, \eg \([10^{-7}, 10^{+1}]\) with \(q = 3 \sim 4\).

\subsubsection{Testing Alternative Trace Estimators}

Besides the Cholesky factorization algorithm, we also repeated the numerical experiments above with stochastic trace estimators, namely, the Hutchinson's algorithm \citep{HUTCHINSON-1990} and stochastic Lanczos quadrature algorithm \citep[Sec. 11.6.1]{GOLUB-2009} to compute the trace of a matrix inverse. These class of randomized methods are attractive due to their scalability to very large matrices, where employing the exact methods could be inefficient, if not infeasible. However, these methods do not compute the exact value of determinant or trace, rather they converge to a solution by Monte-Carlo sampling through iterations. The complexity of Hutchinson method using conjugate gradient is \(\mathcal{O}(\rho n^2 s)\) where \(\rho\) is the density of matrix (\(\rho=1\) for dense matrices) and \(s\) is the number of random vectors for Monte-Carlo sampling. We recall that in our application, the cost of the interpolation scheme is \(q\) times the above-mentioned complexity, \ie
\begin{equation*}
    \mathcal{O} \left( q \rho n^2 s\right).
\end{equation*}
Alternatively, the computational cost of the SLQ method is \(\mathcal{O}( (\rho n^2 + nl) s l )\) where \(l\) is the Lanczos degree, which is the number of Lanczos tri-diagonalization iterations (see details \eg in \cite[Sec. 3]{UBARU-2017}). Thus, the complexity of the interpolation method becomes
    \begin{equation*}
        \mathcal{O} \left( q (\rho n^2 + nl) s l\right).
    \end{equation*}

In both these algorithms, we employed \(s = 30\) random vectors with Rademacher distribution for Monte-Carlo sampling. Also, in SLQ algorithm, we set the Lanczos degree to \(30\).

The results for Hutchinson's algorithm are shown in the fifth to eighth rows, and results for the SLQ algorithm are shown in the ninth to twelfth rows, of \Cref{table:gcv}. Similar to the Cholesky factorization, in both stochastic estimators, the interpolation technique reduces the processing times compared to no interpolation, namely, \(T_{\text{tr}}\) is reduced by two orders of magnitude, and \(T_{\text{tot}}\) by an order of magnitude, while maintaining a reasonable accuracy.

We note that the interpolation with the stochastic methods introduces error due to the uncertainly in the randomized estimation of \(\tau_{-1}\) at the interpolant points \(t_i\). However, the additional error caused by interpolation itself can be less than the error due to aforementioned stochastic estimation. For instance, without interpolation, the SLQ method estimates \(\hat{\vect{y}}\) with a \(1.58 \%\) error, whereas, interpolation with \(q=3\) results in a similar error of \(2.76 \%\) but at a greater than 20-fold reduction in computational cost.


\section{Further Applications} \label{sec:other-app}

We recall that the presented interpolation scheme can be applied to any formulation that consists of the trace or determinant of a  power of the one-parameter affine matrix function \(\tens{A} + t \tens{B}\) where \(\tens{A}\) and \(\tens{B}\) are Hermitian and positive-definite. Often in applications, an algebraic trick (such as in \eqref{eq:trick}) is required to form such an affine matrix function. We here provide two other closely related examples where such affine matrix function can be formulated.


\subsection{Reproducing Kernel Hilbert Space}

Let \(\mathcal{H}_K\) be a reproducing kernel Hilbert space equipped with the reproducing kernel \(K\) that defines the function evaluation \(f(\vect{x}) = \langle f, K(\cdot, \vect{x}) \rangle_{\mathcal{H}_K}\). Consider an infinite-dimensional generalized ridge regression on \(\mathcal{H}_K\) to estimate \(y = f(\vect{x})\) with the given training set \(\{ (\vect{x}_i, y_i) \}_{i=1}^n\) by the minimization problem \citep[Section 5.8.2]{HASTIE-2001}
\begin{equation*}
    \min_{f \in \mathcal{H}_{K}} \, \sum_{i=1}^n \vert y_i - f(\vect{x}_i) \vert^2 + \theta \Vert f \Vert^2_{\mathcal{H}_{K}}.
\end{equation*}
The solution to the above problem has the form \(f(\cdot) = \sum_j \alpha_j K(\cdot, \vect{x}_j)\). For the finite-dimensional formulation, define the kernel matrix \(\tens{K}\) with the components \(K_{ij} \coloneqq K(\vect{x}_i, \vect{x}_j)\), which is symmetric and positive-definite. Let \(\vect{\alpha} \coloneqq [\alpha_1, \dots, \alpha_n]^{\intercal}\) and \(\vect{y} \coloneqq [y_1, \dots, y_n]^{\intercal}\). The minimization problem in finite-dimensional setting becomes
\begin{equation*}
    \min_{\vect{\alpha}} \, \Vert \vect{y} - \tens{K} \vect{\alpha} \Vert_2^2 + \theta \Vert \vect{\alpha} \Vert^2_{\tens{K}},
\end{equation*}
where \(\Vert \vect{\alpha} \Vert^2_{\tens{K}} = \vect{\alpha}^{\intercal} \tens{K} \vect{\alpha}\). The optimal solution to the above problem is
\begin{equation*}
    \hat{\vect{\alpha}} = (\tens{K} + \theta \tens{I})^{-1} \vect{y},
\end{equation*}
and the fitted values on the training points are \(\hat{\vect{y}} = \tens{K} \hat{\vect{\alpha}} =: \tens{S}_{\theta} \vect{y}\) where the smoother matrix \(\tens{S}_{\theta}\) is defined by \(\tens{S}_{\theta} \coloneqq \tens{K}( \tens{K} + \theta \tens{I})^{-1}\).

One may seek the optimal value for \(\theta\) as the minimizer of the GCV function
\begin{equation}
    V(\theta) \coloneqq \frac{\frac{1}{n} \left\| (\tens{I} - \tens{S}_{\theta} ) \vect{y} \right\|_2^2}{\left( \frac{1}{n} \trace( \tens{I} - \tens{S}_{\theta}) \right)^2}. \label{eq:gcv2}
\end{equation}
We recall that the expensive part of computing \eqref{eq:gcv2} is the term \(\trace(\tens{S}_{\theta})\). To apply our interpolation scheme, write \(\tens{S}_{\theta}\) as the \emph{Reinsch} form
\begin{equation*}
    \tens{S}_{\theta} = \theta^{-1} (\tens{K}^{-1} + \theta^{-1} \tens{I})^{-1}.
\end{equation*}
We realize that
\begin{equation*}
    \trace(\tens{S}_{\theta}) = t n (\tau_{-1}(t))^{-1},
\end{equation*}
where \(\t \coloneqq \theta^{-1}\). The proposed interpolation method follows by using \(\tens{A} = \tens{K}^{-1}\), \(\tens{B} = \tens{I}\), and \(\tau_{-1, 0} = n/\trace(\tens{K})\).


\subsection{Kernel-Based GCV for Mixed Models}

Another formulation of kernel-based GCV, for instance by \citet[Equations 9 and 10]{XU-2009}, yields a function of the form
\begin{equation}
    V(h, \phi) = \frac{\frac{1}{n} \Vert (\tens{I} - \tens{H}(h, \phi)) \vect{y} \Vert_2^2}{\left(\frac{1}{n}\trace(\tens{I} - \tens{H}(h, \phi))\right)^2}, \label{eq:gcv3}
\end{equation}
where
\begin{equation}
    \tens{H}(h, \phi) = \tilde{\tens{H}} + (\tens{I} - \tilde{\tens{H}}) \tens{Z} \left(\tens{Z}^{\intercal} (\tens{I} - \tilde{\tens{H}})^{\intercal} (\tens{I} - \tilde{\tens{H}}) \tens{Z} + \gtens{\Sigma} \right)^{-1} \tens{Z}^{\intercal} (\tens{I} - \tilde{\tens{H}})^{\intercal} (\tens{I} - \tilde{\tens{H}}). \label{eq:H}
\end{equation}
In the above, the covariance \(\gtens{\Sigma} = \gtens{\Sigma}(\phi)\) is symmetric and positive-definite, the design matrix of random effects \(\tens{Z}\) has full column-rank, and \(\tilde{\tens{H}} = \tilde{\tens{H}}(h)\) is the smoother matrix when the random effects are absent. Optimal values of the parameters \((h, \phi)\) are sought by minimizing \(V\).

It is possible to represent the term in the denominator of \eqref{eq:gcv3} by the trace of the inverse of a single matrix to be written as \(\tau_{-1}\). To do so, let \(\tens{P} \coloneqq \tens{I} - \tilde{\tens{H}}\) and \(\tens{Y} \coloneqq \tens{P} \tens{Z}\). Using the Woodbury matrix identity and \eqref{eq:H}, we can represent the term inside the trace in \eqref{eq:gcv3} as
\begin{align*}
    \tens{I} - \tens{H}(h, \phi) &= \left( \tens{I} - \tens{Y} (\tens{Y}^{\intercal} \tens{Y} + \gtens{\Sigma})^{-1} \tens{Y}^{\intercal} \right) \tens{P} \\
                                 &= \left( \tens{I} + \tens{Y} \gtens{\Sigma}^{-1} \tens{Y}^{\intercal} \right)^{-1} \tens{P}.
\end{align*}
If \(\tens{P}\) is positive-definite, then let \(\tens{P} = \tens{L} \tens{L}^{\intercal}\) be the Cholesky decomposition of \(\tens{P}\). By using the cyclic property of trace operator, we have
\begin{align}
    \trace(\tens{I} - \tens{H}(h, \phi)) &=\trace\left(\tens{L}^{\intercal} ( \tens{I} + \tens{Y} \gtens{\Sigma}^{-1} \tens{Y}^{\intercal} )^{-1}  \tens{L}\right) \notag \\
        &= \trace\left(( \tens{L}^{-1} \tens{L}^{-\intercal} + \tens{L}^{-1} \tens{Y} \gtens{\Sigma}^{-1} \tens{Y}^{\intercal} \tens{L}^{-\intercal} )^{-1} \right). \label{eq:IH}
\end{align}
Note that both matrices \(\tens{A} \coloneqq \tens{L}^{-1} \tens{L}^{-\intercal}\) and \(\tens{B} \coloneqq \tens{L}^{-1} \tens{Y} \gtens{\Sigma}^{-1} \tens{Y}^{\intercal} \tens{L}^{-\intercal}\) are symmetric and positive-definite since they are in the Gramian form. To compute \eqref{eq:IH}, the presented interpolation method can be applied for instance if \(\gtens{\Sigma}(\phi)\) is linear in its parameter. Such assumption is common, for instance when \(\gtens{\Sigma} = \phi \tens{K}\) where \(\phi\) is variance and \(\tens{K}\) is the correlation matrix. In such a case, the sum of two matrices in \eqref{eq:IH} becomes an affine function of \(t \coloneqq \phi^{-1}\) and \(\trace(\tens{I} - \tilde{\tens{H}}(h, \phi))\) can be written as \(\tau_{-1}(t)\).


\section{Conclusions} \label{sec:conclusion}

In many applications in statistics and machine learning, it is desirable to estimate the determinant and trace of the real powers of a one-parameter family of matrix functions \(\tens{A} + \t \tens{B}\) where the parameter \(\t\) varies and the matrices \(\tens{A}\) and \(\tens{B}\) in the formulation remain unchanged. There exist many efficient techniques to estimate the determinant and trace of implicit matrices (such as inverse of a matrix), however, these methods are geared toward generic matrices. Using those methods, the computation of the determinant and trace of the parametric matrices should be repeated for each parameter value as the matrix is updated. To efficiently perform such computation for a wide range of parameter \(\t\), we presented in this work heuristic methods to interpolate the functions  \(\t \mapsto \log \det(\tens{A} + \t \tens{B})\) and \(\t \mapsto \trace((\tens{A} + \t \tens{B})^{p})\). The interpolation approach is based on sharp bounds for these functions using inequalities for a Schatten-type norm and anti-norm. We proposed two types of interpolation functions, namely, interpolation with a linear combination of orthogonalized inverse-monimial basis functions, and interpolation with rational polynomials, which includes Pad\'{e} approximation and Chebushev rational functions. We demonstrated that both functions can provide highly accurate interpolation over a wide range of \(\t\) using very few interpolation points. The rational polynomials generally provide better results in the neighborhood of the origin of the parameter. In the regions away from the origin, choice of interpolation method is less important; namely we observed e.g., the interpolation with Chebyshev rational functions provide similar results to the orthogonalized inverse-monomials in \eqref{eq:approx} in such cases. All the presented interpolation methods can lead to one to two orders of magnitude savings in processing time in practical applications that require frequent evaluations of \( \log \det(\tens{A} + \t \tens{B})\) or \(\trace((\tens{A} + \t \tens{B})^{p})\). 

For applications where one is interested in values of \(\t \ll \tau_{p, 0}\) (such as in \Cref{sec:example-2} where the matrix was shifted due to being ill-conditioned) interpolation using \eqref{eq:rational} is recommended. One should keep in mind that there exists the possibility that \eqref{eq:rational} can become singular at its poles, but a slight rearrangement of the interpolant points \(\t_i\) can be used to ensure these poles are outside the domain of interest. Although \eqref{eq:rational} provides accurate interpolation for a broad range of \(\t\), for a higher number of interpolation points (e.g 6 or more), relation \eqref{eq:approx} or \eqref{eq:cheb-interp} is preferred. 

In closing, the presented interpolation method can be effectively utilized on large data, particularly with the powerful framework of randomized estimators of trace and log-determinant. A practical application of this method together with stochastic Lanczos quadrature on sparse matrices is given by \citep{AMELI-2022-a} to efficiently train Gaussian process regression. The interested reader may refer to \citep{AMELI-2022-c} where the interpolation scheme can be applied to massive data (\eg \(n \sim 2^{25}\)) using the \texttt{imate} package.

\paragraph{Acknowledgment.}
The authors acknowledge support from the National Science Foundation, award number 1520825, and American Heart Association, award number 18EIA33900046.


 \paragraph{Conflict of interest}
 The authors declare that they have no conflict of interest related to this paper.


\appendix
\renewcommand*{\thesection}{\Alph{section}}

\begin{appendices}

\setcounter{equation}{0}
\renewcommand\theequation{A.\arabic{equation}}


\section{Proofs} \label{sec:pf}

In \Cref{thm:trace-ineq}, we show \eqref{eq:anti-ineq} and \eqref{eq:anti-ineq-rev} for the operator \eqref{eq:schatten} and \(p \in (-\infty, 1) \setminus \{0\}\). The results for \(p = 0\) follows by the continuity condition in \eqref{eq:lim}.

\begin{theorem} \label{thm:trace-ineq}
    Suppose \(p \in (-\infty, 1) \setminus \{0\}\) and let the matrices \(\tens{A}, \tens{B} \in \mathcal{M}_{n, n}(\mathbb{C})\) be Hermitian and positive semi-definite (positive definite if \(p < 0\)). Then
    \begin{subequations}
    \begin{align}
        \| \tens{A} + \tens{B} \|_p &\geq \| \tens{A} \|_p + \| \tens{B}\|_p, \label{eq:trace-ineq} \\
        \| \tens{A} - \tens{B}\|_p &\leq \| \tens{A}\|_p - \| \tens{B}\|_p, \label{eq:trace-ineq-neg}
    \end{align}
    \end{subequations}
        provided that \(\tens{A} \succeq \tens{B}\) (\(\tens{A} \succ \tens{B}\) if \(p < 0\)) for \eqref{eq:trace-ineq-neg} to hold.
    In both \eqref{eq:trace-ineq} and \eqref{eq:trace-ineq-neg}, the equality is achieved if and only if \(\tens{A} = c \tens{B}\) for \(c \in \mathbb{R}_{\geq 0}\) (\(c \in \mathbb{R}_{> 0}\) if \(p < 0\)).
\end{theorem}

We prove \Cref{thm:trace-ineq} as follows.

\begin{definition}[Majorization] \label{def:maj}
    For the \(n\)-tuple \(\vect{x} = (x_1, \dots, x_n) \in \mathbb{R}^n\), we denote by \(\vect{x}^{\downarrow} = (x_1^{\downarrow}, \dots, x_n^{\downarrow})\) the tuple with the same components as \(\vect{x}\), but sorted in decreasing order. Let \(\vect{x}, \vect{y} \in \mathbb{R}^n\). We say \(\vect{x}\) weakly majorizes \(\vect{y}\) from below (submajorizes) and indicate by \(\vect{x} \succ_{\mathrm{w}} \vect{y}\) if and only if
    \begin{equation*}
        \sum_{i=1}^k x_i^{\downarrow} \geq \sum_{i=1}^k y_i^{\downarrow}, \quad \text{for all } k = 1, \dots, n.
    \end{equation*}
    Furthermore, if \(\vect{x} \succ_{\mathrm{w}} \vect{y}\) and \(\sum_{i=1}^n x_i^{\downarrow} = \sum_{i=1}^n y_i^{\downarrow} \), we say \(\vect{x}\) majorizes \(\vect{y}\) and indicate by \(\vect{x} \succ \vect{y}\).
\end{definition}

\begin{proposition} \label{prop:trace-ineq-2}
    Suppose \(p \in (-\infty, 1) \setminus \{0\}\), and let the matrices \(\tens{A}, \tens{B} \in \mathcal{M}_{n, n}(\mathbb{C})\) be Hermitian and positive semi-definite (positive-definite if \(p < 0\)) with the \(n\)-tuple of eigenvalues \(\eig(\tens{A})\) and \(\eig(\tens{B})\), respectively. Then
    \begin{equation}
        \| \tens{A} + \tens{B}\|_p \geq M_p(\eig^{\downarrow}(\tens{A}) + \eig^{\downarrow}(\tens{B})), \label{eq:trace-ineq-2}
    \end{equation}
    where \(M_p\) is the generalized mean defined in \eqref{eq:gen-mean}. The equality in the above holds if and only if \(\eig^{\downarrow}(\tens{A} + \tens{B}) = \eig^{\downarrow}(\tens{A}) + \eig^{\downarrow}(\tens{B})\).
\end{proposition}

\begin{proof}
    We proceed the proof for \(p < 0\) as the case \(p \in (0, 1)\) follows similarly. By Ky Fan eigenvalue inequality for Hermitian matrices \cite[p. 356, Theorems 10.21]{ZHANG-2011},
    \begin{equation}
        \eig(\tens{A} + \tens{B}) \prec \eig^{\downarrow}(\tens{A}) + \eig^{\downarrow}(\tens{B}). \label{eq:kf-1}
    \end{equation}
    Let \(\mathcal{I} \coloneqq [\eig_{\min}, \eig_{\max}]\) where \( \eig_{\min} \coloneqq \min \{ \eig_n^{\downarrow}(\tens{A}), \eig_n^{\downarrow}(\tens{B})\}\) and \(\eig_{\max} \coloneqq \max \{ \eig_1^{\downarrow}(\tens{A}), \eig_1^{\downarrow}(\tens{B})\} \). Since \(\tens{A}\) and \(\tens{B}\) are Hermitian and at least positive semi-definite, we have \(\mathcal{I} \subset \mathbb{R}_{\geq 0}\). Define the convex function \(f(t) \coloneqq t^{p}\) on \(\mathcal{I}\). Applying \eqref{eq:kf-1} to Theorem 5.A.1 of \cite[p. 165]{MARSHALL-2011} yields
    \begin{equation}
        f\left(\eig(\tens{A} + \tens{B})\right) \prec_{\mathrm{w}} f\left(\eig^{\downarrow}(\tens{A}) + \eig^{\downarrow}(\tens{B}) \right). \label{eq:f-week-maj-1}
    \end{equation}
    (We note that for \(p \in (0, 1)\), the function \(f\) is concave and the direction of the above and subsequent inequalities are flipped instead). The above relation in particular, implies
    \begin{equation}
        \sum_{i=1}^n f\left(\eig_i(\tens{A} + \tens{B})\right) \leq \sum_{i=1}^n f\left(\eig_i^{\downarrow}(\tens{A}) + \eig_i^{\downarrow}(\tens{B}) \right). \label{eq:fkf-1}
    \end{equation}
    Raising \eqref{eq:fkf-1} to the power \(\frac{1}{p}\) (which flips the direction of inequality if \(p < 0\)) concludes \eqref{eq:trace-ineq-2}.
    
    Also, the condition \(\eig(\tens{A} + \tens{B})) = \eig^{\downarrow}(\tens{A}) + \eig^{\downarrow}(\tens{B})\) is sufficient for the equality in \eqref{eq:fkf-1}. To show the necessity condition, suppose in contrary that the equality in \eqref{eq:fkf-1} holds. This equality condition together with \eqref{eq:f-week-maj-1} imply
    \begin{equation}
        f\left(\eig(\tens{A} + \tens{B})\right) \prec f\left(\eig^{\downarrow}(\tens{A}) + \eig^{\downarrow}(\tens{B}) \right). \label{eq:fmaj-1}
    \end{equation}
    The above condition can be achieved if and only if either \(f\) is linear \cite[p. 166, Theorem 5.A.1.e]{MARSHALL-2011}, which is not, or if \(\eig(\tens{A} + \tens{B}) = \eig^{\downarrow}(\tens{A}) + \eig^{\downarrow}(\tens{B})\).
\end{proof}

The equality condition in \Cref{prop:trace-ineq-2} is contingent on the following condition.

\begin{lemma} \label{lem:commute}
    Let \(\tens{A}, \tens{B} \in \mathcal{M}_{n, n}(\mathbb{C})\) be Hermitian matrices with the \(n\)-tuple of eigenvalues \(\eig(\tens{A})\) and \(\eig(\tens{B})\), respectively. Then, \( \eig^{\downarrow}(\tens{A} + \tens{B}) = \eig^{\downarrow}(\tens{A}) + \eig^{\downarrow}(\tens{B}) \) only if \(\tens{A}\) and \(\tens{B}\) commute.
\end{lemma}

\begin{proof}
    Using the Golden-Thompson inequality \cite[p. 261, Equation IX.19]{BHATIA-1997} and Von Neumann's trace inequality \citep{MIRSKY-1975} respectively, we have
    \begin{equation*}
        \trace\left(e^{\tens{A} + \tens{B}}\right) \leq \trace\left( e^{\tens{A}} e^{\tens{B}} \right) 
        \leq \sum_{i=1}^n e^{\eig_i^{\downarrow}(\tens{A})} e^{\eig_i^{\downarrow}(\tens{B})} 
        = \sum_{i=1}^n e^{\eig_i^{\downarrow}(\tens{A}) + \eig_i^{\downarrow}(\tens{B})} = \trace\left( e^{\tens{A} + \tens{B}} \right).
    \end{equation*}
    But, the equality in Golden-Thompson inequality holds if and only if \(\tens{A}\) and \(\tens{B}\) commute \citep{PETZ-1994}.
\end{proof}

\begin{remark}
    The equality \eqref{eq:tr-AI} is a special case of \eqref{eq:trace-ineq-2} since \(\tens{A}\) commutes with \(\tens{B} = t \tens{I}\). This also applies to \eqref{eq:det-AI} as it can be obtained from \eqref{eq:tr-AI} at \(p \to 0\).
\end{remark}

We also show superadditivity of the generalized mean function, \(M_p\), for \(p < 1\).


\begin{lemma} \label{lem:concave}
    \(M_p\) for \(p < 1\) is a concave function on \(\mathbb{R}^n_{\geq0}\) (\(\mathbb{R}^n_{> 0}\) if \(p < 0\)).
\end{lemma}

\begin{proof}
    We show the Hessian \(\tens{H}\) of the function \(M_p(\vect{x})\) is negative semi-definite. The component \(H_{ij}\) of the Hessian matrix \(\tens{H}\) is
    \begin{equation*}
        H_{ij} \coloneqq \frac{\partial^2 M_p}{\partial x_i \partial x_j} 
        = \frac{p-1}{n^{\frac{1}{p}}} 
        \left( \sum_{k = 1}^n x_k^p \right)^{\frac{1}{p}-2}
        x_i^{p-1}
        \left( \frac{\delta_{ij}}{x_i} \sum_{k = 1}^n x_k^p - x_j^{p-1} \right),
    \end{equation*}
    where \(\delta_{ij}\) is the Kronecker delta function. The matrix \(\tens{H}\) is negative semi-definite if and only if \(\vect{w}^{\intercal} \tens{H} \vect{w} \leq 0\) for all nonzero vectors \(\vect{w} \coloneqq (w_1,\dots,w_n)\). The latter condition is equivalent to
    \begin{equation*}
        \sum_{i = 1}^n \sum_{j = 1}^n w_i w_j x_i^{p-1} \left( \frac{\delta_{ij}}{x_i} \sum_{k = 1}^n x_k^p - x_j^{p-1} \right) \geq 0,
    \end{equation*}
    which simplifies to
    \begin{equation*}
        \left( \sum_{j = 1}^n w_j x_j^{p-1} \right)^2 \leq \left( \sum_{k = 1}^n x_k^{p} \right) \left( \sum_{i = 1}^n w_i^2 x_i^{p-2} \right). \label{eq:positive-definite}
    \end{equation*}
    The above relation holds by the Cauchy-Schwarz inequality for the product two vectors with the components \(x_j^{\frac{p}{2}}\) and \(w_j x_j^{\frac{p}{2}-1}\). Thus, \(\tens{H}\) is negative semi-definite and it concludes the proof.
\end{proof}

\begin{proposition} \label{prop:ineq}
    \(M_p\) for \(p < 1\) is superadditive on \(\mathbb{R}^n_{\geq0}\) (\(\mathbb{R}^n_{> 0}\) if \(p < 0\)). That is, for \(\vect{x}, \vect{y} \in \mathbb{R}_{>0}^n\),
    \begin{equation}
        M_p(\vect{x} + \vect{y}) \geq M_p(\vect{x}) + M_p(\vect{y}). \label{eq:H+}
    \end{equation}
    The equality in the above holds if and only if \(\vect{x} = c \vect{y}\) where \(c \geq 0\) (\(c > 0\) if \(p < 0\)).
\end{proposition}

\begin{proof}
    Since by \Cref{lem:concave}, the function \(M_p\) is concave, from the Jensen inequality (see \eg \citet[Sec. 3.12]{HARDY-1952}) we have \(M_p( \frac{1}{2}(\vect{x} + \vect{y})) \geq \frac{1}{2} ( M_p(\vect{x}) + M_p(\vect{y}))\), which concludes \eqref{eq:H+}. The Jensen inequality becomes an equality if \(\vect{x} = \vect{y}\). But, since \(M_p(c \vect{x}) = c M_p(\vect{x})\), the equality criterion can be extended to \(\vect{x} = c \vect{y}\).
\end{proof}


We can now prove \Cref{thm:trace-ineq}.

\begin{proof}[Proof of \Cref{thm:trace-ineq}]
    We have \(M_p(\eig^{\downarrow}(\tens{A})) = \| \tens{A} \|_p\) and \(M_p(\eig^{\downarrow}(\tens{B})) = \| \tens{B} \|_p\). Applying \Cref{prop:ineq} to \(M_p(\eig^{\downarrow}(\tens{A}) + \eig^{\downarrow}(\tens{B}))\) and using \Cref{prop:trace-ineq-2} concludes \eqref{eq:trace-ineq}. Also, applying \eqref{eq:trace-ineq} to \(\| \tens{B} + (\tens{A} - \tens{B}) \|_p = \| \tens{A} \|_p\) concludes \eqref{eq:trace-ineq-neg}.

    The equality in \Cref{prop:ineq} holds if and only if \(\eig^{\downarrow}(\tens{A}) = c \eig^{\downarrow}(\tens{B})\) for some positive constant \(c\). Also, by \Cref{lem:commute}, the equality in \Cref{prop:trace-ineq-2} holds if \(\tens{A}\) and \(\tens{B}\) commute, which means \(\tens{A}\) and \(\tens{B}\) have the same eigenspace \citep[p. 50, Theorem 1.3.12]{HORN-1990}. By combining these two conditions, equality in \eqref{eq:trace-ineq} and \eqref{eq:trace-ineq-neg} is achieved when \(\tens{A}\) is a scalar multiple of \(\tens{B}\).
\end{proof}

\end{appendices}


\phantomsection{}
\addcontentsline{toc}{section}{References}



\bibliographystyle{apalike2}

\small
\bibliography{References}


\end{document}